\documentclass[11pt]{elsarticle}
\usepackage{}
\usepackage{amssymb}
\usepackage{xcolor}
\usepackage{mathrsfs}
\usepackage{amscd}
\usepackage{appendix}
\usepackage{geometry}
\usepackage{setspace}
\usepackage{amsfonts}
\usepackage{amsmath}
\usepackage{lineno,hyperref}
\usepackage{tikz}
\usetikzlibrary{matrix,arrows}
\newtheorem{definition}{Definition}[section]
\newtheorem{theorem}[definition]{Theorem} 
\newtheorem{lemma}[definition]{Lemma}     

\newtheorem{remark}[definition]{Remark}

\numberwithin{equation}{section}

\biboptions{numbers,sort&compress} 
\journal{}
\begin{document}
\newgeometry{left=2.5cm,right=2.5cm,top=2.5cm,bottom=2.5cm}

\begin{frontmatter}

\title{Invariant curves of low smooth quasi-periodic reversible mappings }

\author[]{Yan Zhuang}
\ead{yzhuangmath@163.com}

\author[]{Daxiong Piao}
\ead{dxpiao@ouc.edu.cn}

\author[]{Yanmin Niu \corref{cor1}}
\ead{niuyanmin@ouc.edu.cn}
\address{School of Mathematical Sciences,  Ocean University of China, Qingdao 266100, P. R. China}
\cortext[cor1]{Corresponding author}

\begin{abstract}
   In this paper, we obtain the invariant curves of quasi-periodic reversible mappings with finite smoothness. Since the reversible property is difficult to maintain in the process of approximating smooth functions by analytical ones, R\"{u}ssmann's method in \cite{HR} is invalid. Inspired by the recent work of Li, Qi and Yuan in \cite{LJ}, we turn to regard the reversible mapping as the Poincar\'{e} map of a reversible differential equation. By constructing a KAM theorem for a reversible differential equation which is quasi-periodic in time, we obtain the invariant curves of the reversible mapping. Beyond that, we establish some variants of invariant curve theorems for quasi-periodic reversible mappings.
\end{abstract}


\end{frontmatter}

\begin{spacing}{1.25}
\section{Introduction}

In 1962, German mathematician Moser \cite{JM} proved the existence of invariant closed curves on cylinder - Moser's twist theorem. Along with the invariant torus theorems for Hamiltonian systems by former Soviet mathematicians Kolmogorov and Arnold, they form the famous KAM theory. As is well known, the smoothness condition of KAM theory has always been a focus of research. Initially, Kolmogorov and Arnold dealt with analytical systems \cite{CSJM}, while Moser \cite{JM} discussed 333 differentiable area-preserving maps. Is this requirement of 333 smoothness optimal? This issue has attracted the attention of many mathematicians, such as R\"{u}ssmann and Herman. In 1970, R\"{u}ssmann \cite{HR} successfully reduced the smoothness requirement from 333 to only 5. Shortly after, Moser pointed out that R\"{u}ssmann's proof process can actually yield cases with smoothness greater than 3. Interestingly, Herman \cite{MRH1,MRH2} also independently provided proof of this situation. He gave a counterexample to illustrate that for maps that are differentiable less than 3 degrees, there can be no invariant curve. Subsequently, these results were also extended to the quasi-periodic case. Zharnitsky in \cite{ZhV} studied the existence of invariant curves of the following exact symplectic map
 \begin{equation}\label{some-resu}
\begin{cases}
x_1=x+\gamma+y+f(x,y),\\
 y_1=y+g(x,y),
 \end{cases}
\end{equation}
where the perturbations $f,g$ are quasi-periodic in $x$ with the frequency $\omega=(\omega_1,\omega_2,...,\omega_m)$, real analytic, and $\omega_1,\omega_2,...,\omega_m,2\pi\gamma^{-1}$ are sufficiently incommensurable. Moreover, he applied the invariant curve theorem to show the boundedness of solutions of the Fermi-Ulam problem. Recently, Huang \cite{HLL} obtained the existence of invariant curves of planar mappings which are quasi-periodic in the spatial variables, $\mathcal{C}^{p}$ smooth with $p>2n+1$ ($n$ is the number of frequencies) and satisfy the intersection property. The main line of the proofs is similar to that of R\"{u}ssmann \cite{HR1} and the critical step is approximating smooth functions with analytical ones.

The intersection assumption is another issue that KAM theory focuses on. Moser pointed out in \cite{JM} that the existence of invariant curves cannot be guaranteed without making any assumptions. Furthermore, he proposed that the intersection assumption can be replaced by reversibility, and therefore led to the creation of the KAM theory of reversible systems. Sevryuk developed the reversible theory \cite{MBS} for both continuous and discrete systems. Based on the KAM technique, Liu \cite{LB,LB1} proved the existence of invariant curves for real analytic quasi-periodic reversible mappings (\ref{some-resu}). Recently, Zhuang et al. \cite{ZhN} obtained the invariant curves of reversible mappings with all Birkhoff constants being zero by the normal form theory. More works on the invariant curves for reversible mappings may be found in \cite{Li2,LB,LB1,LJ,ZhJ,MBS2} and references therein.

Motivated by above aspects, we tend to study the existence of invariant curves of the standard reversible quasi-periodic mapping (\ref{some-resu}), denoted by $M$,
where mapping $M$ is reversible with respect to the involution $G:(x,y)\rightarrow(-x,y)$, that is, $GMG=M^{-1}$. Moreover, $f,g:\mathbb{R}\times B(r_0)\rightarrow\mathbb{R}$ are $\mathcal{C}^p$ and $\mathcal{C}^{p+m}$, respectively, with $p=2m+1+\mu,0<\mu\ll1$, $B(r_0)=\{|y|\leq r_0:r\in\mathbb{R}\}$, and they are quasi-periodic in $x$ with frequency $\omega=(\omega_1,\omega_2,...,\omega_m)$. That is to say, we are going to obtain the invariant curves of the mappings (\ref{some-resu}) in the smooth cases, other than analytic cases. Such a map is often met when the vector field is quasi-periodic in time and reversible with respect to the involution $G$. The reduction of smoothness for reversible mappings is natural. Moser once gave a conjecture \cite{moser1,moser2} of invariant curves for $\mathcal{C}^4$ reversible map. In fact, for autonomous reversible systems with $\mathcal{C}^4$ functions, P\"{o}schel obtained an existence result of invariant tori in \cite{pochel}. But there is no strict proof whether his results can be generalized to the time-dependent case. Therefore, it is not yet clear whether the twist theorem for $\mathcal{C}^4$ smoothness holds true. Different from the symplectic case in R\"{u}ssmann's work \cite{HR1}, reversibility property is hard to maintain in the process of approximating smooth functions by analytical ones.  Recently Li et al. \cite{LJ} overcame this difficulty by regarding the reversible mapping as the Poincar\'{e} map of a reversible differential equation. They constructed a KAM theorem for higher-dimensional periodic reversible mapping (\ref{some-resu}), where $\gamma$ satisfies Diophantine condition, $f,g$ are $\mathcal{C}^p$ and $\mathcal{C}^{p+d}$ respectively with $p=2d+1+\mu, 0<\mu\ll 1$ and $d$ being the dimension of variables. By this method, we extend the twist theorem to quasi-periodic case, and solve that the existence of invariant curves for non-autonomous reversible system is true for $\mathcal{C}^p$ and $\mathcal{C}^{p+m}$ functions, where $m$ is the dimension of quasi-periodic frequency (see Theorem \ref{theo-map} and \ref{theo-sys}).

In addition, we discuss variant forms of invariant curve theorems. The small twist theorem was designed by Moser \cite{JM} to prove the stability of elliptic fixed points of general type, and it has found many other consequences in stability theory. The result obtained in this paper is useful to simplify the use of the twist theorem in some applications. Based on Moser's work, Ortega \cite{ortega1} obtained a variant of the small twist theorem and also studied the existence of invariant curves of mappings with an average small twist in \cite{ortega}. Later for quasi-periodic analytic reversible mapping, Liu \cite{LB,LB1} established some variants of the invariant curve theorem which are similar to ones in \cite{ortega1} and \cite{ortega}. Subsequently, Huang \cite{HP} extended the results to smooth quasi-periodic mapping with intersection property. Based on our main result, we obtained Theorem \ref{theo-map1} and \ref{theo-map2}.



The rest of the paper is organized as follows. In Section 2, we will list some properties of quasi-periodic functions, and then state the main results. In section 3, we give an approximation Lemma of Jack-Moser-Zhender. The following section is devoted to giving out the key iteration process. In section 5, using Iteration Lemma \ref{iteration}, we give the proof of Theorem \ref{theo-map} and Theorem \ref{theo-sys}. In the last section, we give some variants of our main results.

\section{Quasi-periodic functions and the main results}
\subsection{The space of quasi-periodic functions}

We first define quasi-periodic functions with the frequency $\omega=(\omega_1,\omega_2,...\omega_m)$, here the frequency vector $\omega$ is rational
independent, that is, for all $k=(k_1,k_2,...,k_m)\neq0$, $\langle k,\omega\rangle=\sum_{j=1}^{m}k_j\omega_j\neq 0$.

\begin{definition}\label{defi-shell}
A function $f(t)$ is called a continuous ($\mathcal{C}^p$ or real analytic) quasi-periodic function with the frequency $\omega=(\omega_1,\omega_2,...\omega_m)$, if there is a continuous ($\mathcal{C}^p$ or real analytic) function $F(\theta_1,\theta_2,...,\theta_m)$ which is $2\pi$-periodic in each $\theta_j(1\leq j\leq m)$ such that $$F(\theta)=f(\omega_1t,\omega_2t,...,\omega_mt).$$
This function $F$ is called a shell function of $f$.
\end{definition}

Denote by $Q(\omega)$ the set of real analytic quasi-periodic functions with the frequency $\omega=(\omega_1,\omega_2,...\omega_m)$. From the definition, we know that for given $f(t)\in Q(\omega)$, the corresponding shell function $F:\theta=(\theta_1,\theta_2,...,\theta_m)\in \mathbb{R}^m\rightarrow\mathbb{R}$ has the following Fourier series
$$F(\theta)=\sum_{k\in\mathbb{Z}^m}f_ke^{i\langle k,\theta\rangle},$$
which is a $2\pi$-periodic in each variable, real analytic and bounded in a complex neighbourhood  $\Pi_s^m=\{(\theta_1,\theta_2,...,\theta_m)\in\mathbb{C}^m:|Im\theta_j|\leq s,j=1,2,...,m\}$ of $\mathbb{R}^m$. The function $f$ is obtained from $F$ by replacing $\theta$ by $\omega t$ and has the expansion $$f(t)=\sum_{k\in\mathbb{Z}^m}f_ke^{i\langle k,\omega\rangle t}.$$
\begin{definition}
For $r>0$, let $Q_r(\omega)\subset Q(\omega)$ be the set of real analytic functions $f$ such that the corresponding shell function $F$ are bounded on the subset $\Pi_r^m$, with the supremum norm
$$\|F\|_r=\sup_{\theta\in\Pi_s^m}\left|\sum_{k\in\mathbb{Z}^m}f_ke^{i\langle k,\theta\rangle}\right|=\sup_{\theta\in\Pi_s^m}|F|.$$
Moreover, we define the norm $\|f\|_r=\|F\|_r$.
\end{definition}

The following statement are standard in the theory of quasi-periodic functions and the proof can be found in \cite[chapter 3]{ZhV}.
\begin{lemma}\label{quasi-prop}
The following statement are true:

$(i)$ Let $f(t),g(t)\in Q(\omega)$, then $g(t+f(t))\in Q(\omega)$;

$(ii)$ Suppose $$\left|\langle k,\omega\rangle\right|\geq\frac{c}{|k|^\sigma}, \quad k\in\mathbb{Z}^m\setminus\{0\},$$
for some constants $c>0,\sigma>0$. Let $h(t)\in Q(\omega)$ and $\tau=\alpha t+h(t)(\alpha+h'>0)$, then the inverse relation is given by $t=\alpha^{-1}\tau+h_1(\tau)$ and $h_1\in Q(\frac{\omega}{\alpha})$. In particular, if $\alpha=1$, then $h_1\in Q(\omega)$.
\end{lemma}

\subsection{Invariant curves theorem for quasi-periodic reversible mappings }

 We assume that $f:\mathbb{R}^3\rightarrow\mathbb{R}$ are of class $\mathcal{C}^p$, and define
$$|x|=\max\{|x_1|,|x_2|,|x_3|\},\quad x=(x_1,x_2,x_3)\in\mathbb{R}^3,$$
$$|f|_{\mathbb{R}^3}=\sup_{x\in\mathbb{R}^3}|f(x)|,$$
$$\|f\|_{\mathcal{C}^p}=\sum_{|k|\leq p}\sup_{x\in\mathbb{R}^3}{\left|D^kf(x)\right|},$$
if $p\geq0$ is a integer, and
$$\|f\|_{\mathcal{C}^p}=\sum_{|k|\leq l}\sup_{0<|x- y|<1}\frac{\left|D^kf(x)-D^kf(y)\right|}{|x-y|^s}+\sum_{|k|\leq l}\sup_{x\in\mathbb{R}^3}{\left|D^kf(x)\right|},$$
if $p=l+s,l\geq0$ is an integer, $s\in(0,1)$, where $$D^k=(\frac{\partial}{\partial x_1})^{k_1}\circ(\frac{\partial}{\partial x_2})^{k_2}\circ(\frac{\partial}{\partial x_3})^{k_3},\quad |k|=|k_1|+|k_2|+|k_3|,\quad k=(k_1,k_2,k_3)\in\mathbb{N}^3.$$

Now we are ready to state our main results.
\begin{theorem}\label{theo-map}
Assume that the quasi-periodic mapping $M$ given by (\ref{some-resu}) is reversible with respect to the involution $G:(x,y)\rightarrow(-x,y)$, that is, $M\circ G\circ M=G$.  Given $p=2m+1+\mu$ with $0<\mu\ll 1$, we suppose that $f,g:\mathbb{R}\times B(r_0)\rightarrow\mathbb{R}$ are $\mathcal{C}^p$ and $\mathcal{C}^{p+m}$, respectively. Furthermore, suppose that $\omega_1,\omega_2,...,\omega_m,2\pi\gamma^{-1}$ satisfy the Diophantine condition:
\begin{equation}\label{dioph1}
\left|\langle k,\omega\rangle\frac{\gamma}{2\pi}+l\right|\geq\frac{c_0}{|k|^\sigma},\quad k\in\mathbb{Z}^m\setminus\{0\},l\in\mathbb{Z},
\end{equation}
where $0<c_0<1$, $\sigma=m+\frac{\mu}{100}$.
Then there exists $\varepsilon_0>0$ such that for any $0<\varepsilon<\varepsilon_0$, if
$$\|f\|_{\mathcal{C}^p(\mathbb{R}\times B(r_0))}\leq\varepsilon,\quad\|g\|_{\mathcal{C}^{p+m}(\mathbb{R}\times B(r_0))}\leq\varepsilon,$$
 the mapping $M$ has an invariant curve.
\end{theorem}
\begin{theorem}\label{theo-sys}
Consider a system of non-autonomous differential equations
\begin{equation}\label{main-sys}A:
\begin{cases}
x'=\gamma+y+f(x,y,t),\\
y'=g(x,y,t)
\end{cases}
(x,y,t)\in D:=\mathbb{R}\times B(r_0)\times\mathbb{T}.
\end{equation}
Suppose that the system $A$ is reversible with respect to the involution $G:(x,y,t)\rightarrow(x,y,-t)$, that is, for any $(x,y,t)\in D$,
$$f(-x,y,-t)=f(x,y,t),\quad g(-x,y,-t)=-g(x,y,t).$$
Given $p=2m+1+\mu$ with $0<\mu\ll 1$, we assume that $f,g:D\rightarrow\mathbb{R}$ are $\mathcal{C}^p$ and $\mathcal{C}^{p+m}$, respectively. Moreover, assume that $\omega$ and $\gamma^{-1}$ satisfy the Diophantine condition:
\begin{equation}\label{dioph2}
\left|\langle k,\omega\rangle\gamma+l\right|\geq\frac{c_0}{|k|^\sigma},\quad k\in\mathbb{Z}^m\setminus\{0\},l\in\mathbb{Z},
\end{equation}
where $0<c_0<1$, $\sigma=m+\frac{\mu}{100}$.
Then there exists $\varepsilon_0>0$ such that for any $0<\varepsilon<\varepsilon_0$, if
$$\|f\|_{\mathcal{C}^p(D)}\leq\varepsilon,\quad\|g\|_{\mathcal{C}^{p+m}(D)}\leq\varepsilon,$$
 the system  $A$ has an invariant curve.
\end{theorem}
\begin{remark}\label{remark}
The first equation of system (\ref{main-sys}) can be replaced by $\gamma+h(y)+f(x,y,t)$ with $h'(y)> 0$.
\end{remark}

\subsection{The small twist theorem}

In this section, we are ready to give a useful small twist theorem which is a variant of the invariant curve theorem (Theorem \ref{theo-map}) for the quasi-periodic reversible mapping $M$.

In many applications, one may meet the so called small twist mappings. We consider the following reversible mappings $M_1:$
\begin{equation}
\begin{cases}
x_1=x+\gamma+\delta y+\delta f(x,y;\delta),\\
y_1=y+\delta g(x,y;\delta),
\end{cases}
\end{equation}
where functions $f$ and $g$ are quasi-periodic in $x$ with the frequency $\omega=(\omega_1,\omega_2,...,\omega_m)$.  $ f(x,y,0)= g(x,y,0)=0$, $0<\delta<1$ is a small parameter.
\begin{theorem}\label{theo-map1}
Assume that the quasi-periodic mapping $M_1$ is reversible with respect to the involution $G:(x,y)\rightarrow(-x,y)$, that is, $M_1\circ G\circ M_1=G$.  Given $p=2m+1+\mu$ with $0<\mu\ll 1$, we suppose that $f,g:\mathbb{R}\times B(r_0)\rightarrow\mathbb{R}$ are $\mathcal{C}^p$ and $\mathcal{C}^{p+m}$, respectively. Furthermore, suppose that $\omega_1,\omega_2,...,\omega_m,2\gamma^{-1}\pi$ satisfy the Diophantine condition (\ref{dioph1}).
Then there are two positive numbers $\Delta_0$ and $\varepsilon_0$ such that if $0<\delta<\Delta_0$
and \begin{equation}\label{perturbation}
\|f\|_{\mathcal{C}^p(\mathbb{R}\times B(r_0))}\leq\varepsilon_0,\quad\|g\|_{\mathcal{C}^{p+m}(\mathbb{R}\times B(r_0))}\leq\varepsilon_0,
\end{equation}
 the mapping $M_1$ has an invariant curve.
\end{theorem}
\begin{remark}
Theorem \ref{theo-map1} is the so called small twist theorem. The proof of this theorem when $f, g$ are real analytic can be found in \cite{MBS}, as far as we all know, there is no proof if these functions are $\mathcal{C}^p$ smooth. It is not a direct consequence of Theorem \ref{theo-map}, but one can use the same procedure in the proof of Theorem \ref{theo-map} to prove it. Thus we omit the proof of Theorem \ref{theo-map1}.
\end{remark}
\begin{remark}\label{zhu}
The above results is also true for the following quasi-periodic reversible mapping $M_2$:
\begin{equation}
\begin{cases}
x_1=x+\gamma+\delta h(y)+\delta f(x,y;\delta),\\
y_1=y+\delta g(x,y;\delta),
\end{cases}
\end{equation}
with $h'(y)>0$. That is to say, assume the conditions of Theorem  \ref{theo-map1} hold, then there exist two positive numbers $\Delta_0$ and $\varepsilon_0$, such that if $0<\delta<\Delta_0$ and
\begin{equation}\label{perturbation1}
\|f\|_{\mathcal{C}^p(\mathbb{R}\times B(r_0))}\leq\varepsilon_0,\quad\|g\|_{\mathcal{C}^{p+m}(\mathbb{R}\times B(r_0))}\leq\varepsilon_0,
\end{equation}
 the mapping $M_2$ has an invariant curve.
\end{remark}

In the following, we are going to give a variant of the small twist theorem. Consider a one-parameter family of mappings $\{M_\delta\}_{\delta\in[0,1]}$ with $M_\delta: \mathbb{R}\times B(r_0)\rightarrow \mathbb{R}\times\mathbb{R}$. $M_\delta$ can be described in the form:
\begin{equation}
\begin{cases}
x_1=x+\gamma+\delta l_1(x,y)+\delta f(x,y,\delta),\\
y_1=y+\delta l_2(x,y)+\delta g(x,y,\delta),
\end{cases}
\end{equation}
where $ f(x,y,0)= g(x,y,0)=0$, $l_1,l_2,f,g$ are quasi-periodic in $x$ with the frequency $\omega=(\omega_1,\omega_2,...,\omega_m)$. $0<\delta<1$ is a small parameter.

Since $l_1$ and  $l_2$ depend on the angle variable $x$, it seems that one cannot use this result directly to infer the existence of invariant curves of mapping  $M_\delta$. However, if $f,g$  satisfy some further conditions and under series of attempts, we can construct a change of variables such that the original mapping $M_\delta$ is transformed into a new one, which has the same form as $M_2$. At the same time, the new mapping meets all conditions of Remark \ref{zhu}. More precisely, we will prove the following results.
\begin{theorem}\label{theo-map2}
Assume that the quasi-periodic mapping $M_\delta$ is reversible with respect to the involution $G:(x,y)\rightarrow(-x,y)$, that is, $M_\delta\circ G\circ M_\delta=G$.  Given $p=2m+1+\mu$ with $0<\mu\ll 1$, we suppose that $f:\mathbb{R}\times B(r_0)\rightarrow\mathbb{R}$ are $\mathcal{C}^{p}$ and $l_1,l_2,g:\mathbb{R}\times B(r_0)\rightarrow\mathbb{R}$ are $\mathcal{C}^{p+m}$. Furthermore, Suppose that $\omega_1,\omega_2,...,\omega_m,2\gamma^{-1}\pi$ satisfy the Diophantine condition (\ref{dioph1}).
In the meanwhile,
\begin{equation}\label{average1}
\lim\limits_{T\rightarrow\infty}\frac{1}{T}\int_0^T\frac{\partial l_1}{\partial y}(x,y)dx> 0.
\end{equation}
Then there are two positive numbers $\bar{\Delta}_0$ and $\varepsilon$ such that if $0<\delta<\bar{\Delta}_0$ and
\begin{equation}\label{perturbation2}
\|f\|_{\mathcal{C}^p(\mathbb{R}\times B(r_0))}\leq\varepsilon,\quad\|g\|_{\mathcal{C}^{p+m}(\mathbb{R}\times B(r_0))}\leq\varepsilon,
\end{equation}
 the mapping $M_\delta$ has an invariant curve.
\end{theorem}
\section{Approximation Lemma}

In this section, we will give a well known and fundamental approximation results, which is used in the iterative process.

First we define the kernel function
$$K(w)=\frac{1}{(2\pi)^l}\int_{\mathbb{R}^l}\hat{K}(\xi)e^{i\langle w,\xi\rangle}d\xi,\quad w\in{\mathbb{C}^l},$$
where $\hat{K}$ is a $\mathcal{C}^\infty$ function with compact support, contained in the ball $|\xi|\leq a$ with a constant $a>0$, that satisfies
\begin{equation*}
\partial^k\hat{K}(0)=
\begin{cases}
1,\quad if~k=0,\\
0,\quad if ~k\neq 0.
\end{cases}
\end{equation*}
\begin{lemma}(Jackson-Moser-Zehnder)\label{JMZ} Let $f(w)\in\mathcal{C}^p(\mathbb{R}^l)$ for some $p>0$ with finite $\mathcal{C}^p$ norm over $\mathbb{R}^l$. For any $\delta>0$, define
\begin{equation}
(S_\delta f)(w)=\delta^{-l}\int_{\mathbb{R}^l}K\Big(\frac{w-w^*}{\delta}\Big)f(w^*)dw^*.
\end{equation}
Then there exists a constant $c\geq1$ depending only on $p$ and $l$ such that for any $\delta>0$, the function $(S_\delta f)(w)$ is real analytic on $\mathbb{C}^l$, and for any $k\in\mathbb{N}^l$ with $|k|\leq p$, one has
$$\sup_{w\in\Pi_\delta^l}\left|\partial^k(S_\delta f)(w)-\sum_{|\lambda|\leq p-|k|}\frac{\partial^{\lambda+k}(S_\delta f)(Rew)}{\lambda!}(iImw)^\lambda\right|\leq C\|f\|_{\mathcal{C}^p}\delta^{p-|k|},$$
and for all $0<\delta<\delta'$,
$$\sup_{w\in\Pi_\delta^l}\left|\partial^k(S_{\delta'} f)(w)-\partial^k(S_\delta f)(w)\right|\leq C\|f\|_{\mathcal{C}^p}\delta'^{p-|k|}.$$
Moreover, the H$\ddot{o}$lder norms of $S_\delta f$ satisfy, for all $0\leq p'\leq p\leq p''$,
$$\|S_\delta f-f\|_{\mathcal{C}^{p'}}\leq C\|f\|_{\mathcal{C}^p}\delta^{p-p'},\quad \|S_\delta f\|_{\mathcal{C}^{p''}}\leq C\|f\|_{\mathcal{C}^p}\delta^{p-p''}.$$
Finally, the function $S_\delta f$ preserves periodicity, that is, if $f$ is $T$-periodic in any of its variables $w_j(1\leq j\leq l)$, so is $S_\delta f$.
\end{lemma}

See \cite{LG,DS,EZ} for the proof of Lemma \ref{JMZ}.
\begin{lemma}\label{analy-approxi}
Let $l=3$ in Lemma \ref{JMZ}, assume that $f(x,y,t)\in\mathcal{C}^p$ is a quasi-periodic function with the frequency $\omega=(\omega_1,\omega_2,...,\omega_m)$.
Then for any $\delta>0$, there exists a holomorphic function $S_\delta f:\mathbb{C}\times\mathbb{C}\times\mathbb{T}\rightarrow\mathbb{C}, S_\delta f(\mathbb{R}^3)\subseteq\mathbb{R}$ such that the following inequalities
\begin{equation}\label{inequa}
\begin{cases}
|S_\delta f|_{D_\delta}\leq C\|f\|_{\mathbb{R}^3},\\
|S_\delta f-f|_{\mathbb{R}^3}\leq C\|f\|_{\mathcal{C}^p}\delta^p,\\
|S_\delta f-S_{\delta'} f|_{D_\delta}\leq C\|f\|_{\mathcal{C}^p}\delta'^p
\end{cases}
\end{equation}
hold for $0<\delta<\delta'$, where $D_\delta=\mathbb{R}_\delta\times\mathbb{R}_\delta\times\mathbb{T}_\delta, \mathbb{R}_\delta=\{x\in\mathbb{C}:|Imx|<\delta\},\mathbb{T}_\delta=\{t\in(\mathbb{C}/{2\pi\mathbb{Z}}):|Im t|<\delta\}$.  Moreover, for $(x,y,t)\in\mathbb{R}\times\mathbb{R}\times\mathbb{T}$,

$\bullet$ if $f(-x,y,-t)=f(x,y,t)$, then
\begin{equation}\label{even}S_\delta f(-x,y,-t)=S_\delta f(x,y,t);
\end{equation}

$\bullet$ if $f(-x,y,-t)=-f(x,y,t)$, then
\begin{equation}\label{odd}S_\delta f(-x,y,-t)=-S_\delta f(x,y,t);
\end{equation}
\end{lemma}
{\bf Proof} ~Since that $f(x,y,t)$ is quasi-periodic in $\theta$ with the frequency $\omega=(\omega_1,\omega_2,...,\omega_m)$, from Definition \ref{defi-shell}, there exists a corresponding shell function $$F(\theta,y,t)=F(\theta_1,\theta_2,...,\theta_m,y,t),\quad \theta=(\theta_1,\theta_2,...,\theta_m),$$
which is $2\pi$-periodic in $\theta_j(j=1,2,...,m)$ such that $f(x,y,t)=F(\omega_1 x,\omega_2x,...,\omega_mx,y,t)$.
Since $f(x,y,t)\in\mathcal{C}^p(\mathbb{R}^3)$, it yields $F\in\mathcal{C}^p(\mathbb{R}^{m+2})$. From Lemma \ref{JMZ}, there exists a real analytic function $S_\delta F(\theta_1,\theta_2,...,\theta_m,y,t)$ which is $2\pi$-periodic in $\theta_j(j=1,2,...,m)$ and $t$. It means that there is a real analytic and quasi-periodic function $S_\delta f(x,y,t)$ with the frequency $\omega=(\omega_1,\omega_2,...,\omega_m)$. The proof of inequalities of (\ref{inequa}) is similar to Lemma 2.11 in \cite{HLL}, here we omit it.

Let $w=(w_1,w_2,w_3)\in\mathbb{R}^3$, choose the kernel function $\hat{K}(w)$ such that
$$\hat{K}(-w_1,w_2,-w_3)=\hat{K}(w_1,w_2,w_3), \quad\forall(w_1,w_2,w_3)\in\mathbb{R}^3.$$
By the definition of $K$, it is easy to obtain that $$K(-w_1,w_2,-w_3)=K(w_1,w_2,w_3).$$
Since $$(S_\delta f)(x,y,t)=\delta^{-3}\int_{\mathbb{R}^3}K\Big(\delta^{-1}(x-\tilde{x},y-\tilde{y},t-\tilde{t})\Big)f(\tilde{x},\tilde{y},\tilde{t})d\tilde{x}d\tilde{y}d\tilde{t},$$
\begin{equation*}
\begin{aligned}
(S_\delta f)(-x,y,-t)&=\delta^{-3}\int_{\mathbb{R}^3}K\Big(\delta^{-1}(-x-\tilde{x},y-\tilde{y},-t-\tilde{t})\Big)f(\tilde{x},\tilde{y},\tilde{t})d\tilde{x}d\tilde{y}d\tilde{t}\\
&=\delta^{-3}\int_{\mathbb{R}^3}K\Big(\delta^{-1}(x+\tilde{x},y-\tilde{y},t+\tilde{t})\Big)f(\tilde{x},\tilde{y},\tilde{t})d\tilde{x}d\tilde{y}d\tilde{t}\\
&=\delta^{-3}\int_{\mathbb{R}^3}K\Big(\delta^{-1}(x-\tilde{x}^*,y-\tilde{y}^*,t-\tilde{t}^*)\Big)f(-\tilde{x}^*,\tilde{y}^*,-\tilde{t}^*)d\tilde{x}^*d\tilde{y}^*d\tilde{t}^*\\
&=\delta^{-3}\int_{\mathbb{R}^3}K\Big(\delta^{-1}(x-\tilde{x}^*,y-\tilde{y}^*,t-\tilde{t}^*)\Big)f(\tilde{x}^*,\tilde{y}^*,\tilde{t}^*)d\tilde{x}^*d\tilde{y}^*d\tilde{t}^*\\
&=(S_\delta f)(x,y,t).
\end{aligned}
\end{equation*}
Hence, (\ref{even}) holds. According to the same process, we can also obtain that (\ref{odd}) holds. Thus we complete the proof of this Lemma.$\hfill\square$

Consider a $\mathbb{R}$-value function $f:D\rightarrow\mathbb{R}$. By Whitney's extension theorem, we can find a $\mathbb{R}$-value function $\tilde{f}:\mathbb{R}\times\mathbb{R}\times\mathbb{T}\rightarrow\mathbb{R}$, such that $\tilde{f}|_{D}=f$ and
$$\|\tilde{f}\|_{\mathcal{C}^{k}(\mathbb{R}\times\mathbb{R}\times\mathbb{T})}\leq C\|f\|_{\mathcal{C}^{k}(D)},\quad \forall k\in\mathbb{N},k\leq p.$$

Fix a sequence of fast decreasing $s_n\downarrow 0$, $n\in\mathbb{Z}$ and $s_0\leq\frac{1}{2}$. Let $\delta=s_n$, from Lemma \ref{analy-approxi},  there exists a sequence $\{f^n(w)\}_{n=0}^{\infty}, w=(x,y,t)$, where $$f^n(w)=(S_{s_n} \tilde{f})(w).$$
Moreover, $f^n(w)$ obey the following properties:

$(i):$ $f^n(w)(n=0,1,...)$ are real analytic and quasi-periodic functions with the frequency $\omega=(\omega_1,\omega_2,...,\omega_m)$ on the complex domain $D_{s_n}=\mathbb{R}_{s_n}\times\mathbb{R}_{s_n}\times\mathbb{T}_{s_n}$. In the sequel, denote $D_{s_n}$ as $D_n$ for short.

$(ii):$ $f^n(w)(n=0,1,...)$ satisfy the following inequalities
$$\sup_{w\in D_0}|f^0(w)|\leq C\|f\|_{\mathcal{C}^p(D)},$$
$$\sup_{w\in D_n}|f^n(w)-f(w)|\leq C\|f\|_{\mathcal{C}^p(D)}s_n^p,$$
$$\sup_{w\in D_{n+1}}|f^{n+1}(w)-f^n(w)|\leq C\|f\|_{\mathcal{C}^p(D)}s_n^p,$$
where constant $C$ depend on only $m$ and $p$.

Let $$f_0(w)=f^0(w), f_{n+1}(w)=f^{n+1}(w)-f^{n}(w),$$ we have
\begin{equation*}
f(w)=f^0(w)+\sum_{n=0}^{\infty}\Big(f^{n+1}(w)-f^{n}(w)\Big)=\sum_{n=0}^\infty f_{n}(w),\quad w\in D.
\end{equation*}
By Lemma \ref{analy-approxi}, if $f(-x,y,-t)=f(x,y,t)$, then
\begin{equation*}
f_{n}(-x,y,-t)=f_{n}(x,y,t);
\end{equation*}
 if $f(-x,y,-t)=-f(x,y,t)$, then
\begin{equation*}
f_{n}(-x,y,-t)=-f_{n}(x,y,t).
\end{equation*}

\section{The iteration process}

In this section, we are devoted to present an iteration process leading to the proof of Theorem \ref{theo-map} and \ref{theo-sys}. Before we describe the iteration process, we set up some constants and notations.
\subsection{Constants and Notations}
$(i):$ Given constant $0<\mu\ll 1$, denote $\sigma=m+\frac{\mu}{100}$, $\tilde{\mu}=\frac{\mu}{100(2\tau+1+\mu)}$ and $p=2m+1+\mu$;

$(ii):$
\begin{equation*}
\begin{cases}
\varepsilon_0=\varepsilon,\quad \varepsilon_n=\varepsilon^{(1+\tilde{\mu})^n},\\
s_n=\varepsilon_n^{\frac{1}{p}},\quad r_n=s_n^{m+1+\frac{\mu}{10}},\tau_n=\varepsilon^{-(1+\tilde{\mu})^{n-1}\cdot(1+\frac{m}{p})\tilde{\mu}},\\
s_n^j=s_n-\frac{j}{100p}(s_n-s_{n+1}),\quad r_n^j=r_n-\frac{j}{100p}(r_n-r_{n+1});
\end{cases}
\end{equation*}

$(iii):$ $\mathbb{B}_{\mathbb{C}}(r)=\{y\in\mathbb{C}:|y|\leq r,r\geq0\}$;

$(iv):$ $D(s,r)=\mathbb{R}_s\times\mathbb{B}_{\mathbb{C}}(r)\times\mathbb{T}_s$.

We use a norm definition $\|f\|_{s,r}=\sup\limits_{w=(x,y,t)\in D(s,r)}|f(w)|$.

Now we are going back to the system (\ref{main-sys}). By the discussion of previous section, the system (\ref{main-sys}) can turn into
\begin{equation}
\begin{cases}
x'=\omega+y+\sum\limits_{n=0}^\infty f_n(x,y,t),\\
y'=\sum\limits_{n=0}^\infty g_n(x,y,t),
\end{cases}
\end{equation}
where functions
$f_n,g_n:D_n\rightarrow\mathbb{C}$ are real function and quasi-periodic with the frequency $\omega=(\omega_1,\omega_2,...,\omega_m)$, satisfying
\begin{equation}
\|f_n\|_{s_n,r_n}\leq C\varepsilon\cdot\varepsilon_{n-1},\quad \|g_n\|_{s_n,r_n}\leq C\varepsilon\cdot\varepsilon_{n-1}s_{n-1}^{m}.
\end{equation}
Moreover, 
\begin{equation}\label{even-odd}
f_n(-x,y,-t)=f_n(x,y,t),\quad g_n(-x,y,-t)=-g_n(x,y,t),\quad \forall (x,y,t)\in D(s_n,r_n).
\end{equation}

\subsection{Iteration Lemma}
\begin{lemma}\label{iteration}
Assume that we have $n$ coordinate changes $\Delta\Phi_1,\Delta\Phi_2,...,\Delta\Phi_{n}$.
$\Delta\Phi_j:D(s_j,r_j)\rightarrow D(s_{j-1},r_{j-1})(j=1,2,...,n)$ has the form
\begin{equation}
x=\xi+u_j(\xi,\eta,t),\quad y=\eta+v_j(\xi,\eta,t),
\end{equation}
where $u_j,v_j$ are real analytic and quasi-periodic in $\xi$ with the frequency $\omega$ such that $\Phi_j$ is $G-$invariant with respect to $G(\xi,\eta,t)=(-\xi,\eta,-t)$. Moreover,
\begin{align}
&\|u_j\|_{s_j,r_j}\leq C\varepsilon\cdot\tau_{n-1}\varepsilon_{n-1}s_{n-1}^{-(m+\frac{11\mu}{100})},\label{estimate1}\\
& \|v_j\|_{s_j,r_j}\leq C\varepsilon\cdot\tau_{n-1}\varepsilon_{n-1}s_{n-1}^{-\frac{\mu}{100}}.\label{estimate1-1}
\end{align}
Under the transformation $\Phi_n=\Delta\Phi_1\circ\Delta\Phi_2\circ\cdot\cdot\cdot\circ\Delta\Phi_n$,  the system
\begin{equation}\label{initial-sys}
A^{n-1}:
\begin{cases}
x'=\gamma+y+\sum\limits_{j=0}^{n-1}f_j(x,y,t),\\
y'=\sum\limits_{j=0}^{n-1}g_j(x,y,t)
\end{cases}
\end{equation}
 is changed into
\begin{equation}\label{modified-sys1}
(\partial \Phi_n)^{-1}A^{n-1}\circ\Phi_n:
\begin{cases}
\xi'=\gamma+\eta+\bar{f}_n(\xi,\eta,t),\\
\eta'=\bar{g}_n(\xi,\eta,t),
\end{cases}
\end{equation}
where $\bar{f}_n,\bar{g}_n$ have the following properties:

$(i):$ $\bar{f}_n,\bar{g}_n$ are real analytic and quasi-periodic in $\xi$ with the frequency $\omega$;

$(ii):$ $\bar{f}_n,\bar{g}_n$ has the estimates
\begin{equation}\label{estimate2}
\|\bar{f}_n\|_{s_n,r_n}\leq C\varepsilon\cdot\varepsilon_{n-1},\quad \|\bar{g}_n\|_{s_n,r_n}\leq C\varepsilon\cdot\varepsilon_{n-1}s_{n-1}^{m};
\end{equation}

$(iii):$
\begin{equation}\label{estimate3}
\bar{f}_n(-\xi,\eta,-t)=\bar{f}_n(\xi,\eta,t),\quad\bar{g}_n(-\xi,\eta,-t)=-\bar{g}_n(\xi,\eta,t).
\end{equation}

Then there exists a coordinate change $\Delta\Phi_{n+1}:D(s_{n+1},r_{n+1})\rightarrow D(s_n,r_n)$, which has the form
$$x=\xi+u_{n+1}(\xi,\eta,t),\quad y=\eta+v_{n+1}(\xi,\eta,t),$$
where $u_{n+1},v_{n+1}$ are real analytic and quasi-periodic in $\xi$ with the frequency $\omega$ such that $\Delta\Phi_{n+1}$ is $G-$invariant with respect to $G(\xi,\eta,t)=(-\xi,\eta,-t)$. Moreover,
\begin{align}
&\|u_{n+1}\|_{s_{n+1},r_{n+1}}\leq C\varepsilon\cdot\tau_{n}\varepsilon_{n}s_{n}^{-(m+\frac{11\mu}{100})},\label{estimate4}\\
&\|v_{n+1}\|_{s_{n+1},r_{n+1}}\leq C\varepsilon\cdot\tau_{n}\varepsilon_{n}s_{n}^{-\frac{\mu}{100}}.\label{estimate4-1}
\end{align}
Under the transformation $\Delta\Phi_{n+1}$, the modified system of (\ref{modified-sys1}):
\begin{equation}\label{modified-sys2}
A_{n}:
\left(
\begin{array}{c}
x'\\
y'
\end{array}
\right)
=
\left(
\begin{array}{ccc}
\gamma+y+\bar{f}_n(x,y,t)\\
\bar{g}_n(x,y,t)
\end{array}
\right)
+
(\partial \Phi_n)^{-1}(f_n,g_n)\circ\Phi_n
\end{equation}
 is changed into
\begin{equation}\label{modified-sys3}
(\partial \Delta\Phi_{n+1})^{-1}A_n\circ\Delta\Phi_{n+1}:
\begin{cases}
\xi'=\gamma+\eta+\bar{f}_{n+1}(\xi,\eta,t),\\
\eta'=\bar{g}_{n+1}(\xi,\eta,t),
\end{cases}
\end{equation}
where $\bar{f}_n,\bar{g}_n$ obey the conditions $(i)-(iii)$ by replacing $n$ by $n+1$. Then under the transformation $\Phi_{n+1}=\Phi_n\circ\Delta\Phi_{n+1}$, the system $A^n$ is changed into
\begin{equation}\label{modified-sys4}
(\partial \Phi_{n+1})^{-1}A^{n}\circ\Phi_{n+1}=(\partial \Delta\Phi_{n+1})^{-1}A_n\circ\Delta\Phi_{n+1}.
\end{equation}
\end{lemma}
\subsection{The proof of Iteration Lemma}

In this section, we will first guarantee the existence of the change $\Delta\Phi_{n+1}$ by determining the unknown functions $u_{n+1}(\xi,\eta,t), v_{n+1}(\xi,\eta,t)$, then estimate $u_{n+1}(\xi,\eta,t),v_{n+1}(\xi,\eta,t)$ and new disturbances.

$(\mathrm{I})$~First of all, define
\begin{equation*}
(\hat{f}_n,\hat{g}_n)=(\bar{f}_n,\bar{g}_n)+(\partial \Phi_n)^{-1}(f_n,g_n)\circ\Phi_n,
\end{equation*}
then (\ref{modified-sys2}) can be written as
\begin{equation}\label{modified-sys5}
A_{n}:
(x',y')=\big(\gamma+y+\hat{f}_n(x,y,t),\hat{g}_n(x,y,t)\big).
\end{equation}
By the assumptions and the condition $(i)$ in the iteration lemma, we derive that

$(1):$ $\hat{f}_n,\hat{g}_n$ are real analytic and quasi-periodic in $\xi$ with the frequency $\omega$;

By the assumptions and the condition $(ii)$ in the iteration lemma, it is easy to obtain that
\begin{equation}\label{change-esti1}
\|\Delta\Phi_{n}-Id\|_{s_n,r_n}\leq C\varepsilon\cdot\tau_{n-1}\varepsilon_{n-1}s_{n-1}^{-(m+\frac{11\mu}{100})}.
\end{equation}
Observe that
\begin{equation}\label{change-esti2}
\begin{aligned}
\|\partial \Phi_{n}\|_{s_n,r_n}=&\|(\partial\Delta\Phi_1\circ\Delta\Phi_2\circ\cdot\cdot\cdot\circ\Delta\Phi_n)(\partial\Delta\Phi_2\circ\cdot\cdot\cdot\circ\Delta\Phi_n)\cdots (\partial\Delta\Phi_n)\|_{s_n,r_n}\\
\leq&\prod\limits_{j=1}^{n}\Big(1+C\varepsilon\cdot\tau_{j-1}\varepsilon_{j-1}s_{j-1}^{-(m+\frac{11\mu}{100})}\max\{s_{j}^{-1},r_{j}^{-1}\}\Big)\\
\leq&\prod\limits_{j=1}^{n}(1+C\varepsilon_{j-1}^{\frac{39\mu}{50p}})\leq C.
\end{aligned}
\end{equation}
Combing with the definition of $\hat{f}_n,\hat{g}_n$, it yields that

$(2):$ $\hat{f}_n,\hat{g}_n$ has the estimates
\begin{equation}\label{estimate5}
\|\hat{f}_n\|_{s_n,r_n}\leq C\varepsilon\cdot\varepsilon_{n-1},\quad \|\hat{g}_n\|_{s_j,r_j}\leq C\varepsilon\cdot\varepsilon_{n-1}s_{n-1}^{m};
\end{equation}

By the conditions $(iii)$ and (\ref{even-odd}), we have

$(3):$
\begin{equation}\label{estimate6}
\hat{f}_n(-x,y,-t)=\hat{f}_n(x,y,t),\quad\hat{g}_n(-x,y,-t)=-\hat{g}_n(x,y,t).
\end{equation}

Assume that the change $\Delta\Phi_{n+1}=\Phi$ has the form:
\begin{equation}\label{change1}
\begin{cases}
x=\xi+u(\xi,\eta,t),\\
y=\eta+v(\xi,\eta,t)
\end{cases}
\end{equation}
and its inverse $\Psi=\Phi^{-1}$ is of the form:
\begin{equation}\label{change2}
\begin{cases}
\xi=x+u^*(x,y,t),\\
\eta=y+v^*(x,y,t)
\end{cases}
\end{equation}
where $u,v,u^*,v^*$ are determined later.
Differentiating (\ref{change2}), we have that
\begin{equation}\label{change3}
\xi'=x'+\partial _x u^*\cdot x'+\partial _y u^*\cdot y'+\partial _t u^*,\quad
\eta'=y'+\partial _x v^*\cdot x'+\partial _y v^*\cdot y'+\partial _t v^*.
\end{equation}
Inserting (\ref{modified-sys5}) into (\ref{change3}), we get
\begin{equation}\label{change4}
\begin{aligned}
\xi'=&\gamma+\eta\\
+&\gamma\partial _x u^*+\partial _t u^*+\hat{f}_n(x,y,t)-v^*(x,y,t)\\
+&\partial _y u^*\cdot\hat{g}_n(x,y,t)+\partial _x u^*\cdot\hat{f}_n(x,y,t)\\
+&\partial _x u^*\cdot y
\end{aligned}
\end{equation}
and
\begin{equation}\label{change5}
\begin{aligned}
\eta'=&\gamma\partial _x v^*+\partial _t v^*+\hat{g}_n(x,y,t)\\
+&\partial _y v^*\cdot\hat{g}_n(x,y,t)+\partial _x v^*\cdot\hat{f}_n(x,y,t)\\
+&\partial _xv^*\cdot y.
\end{aligned}
\end{equation}

Consider homological equations:
\begin{equation}\label{homo-equa}
\begin{cases}
\gamma\partial _x u^*+\partial _t u^*+\hat{f}_n(x,y,t)-v^*(x,y,t)=0,\\
\gamma\partial _x v^*+\partial _t v^*+\hat{g}_n(x,y,t)=0.
\end{cases}
\end{equation}

We first solve the function $v^*$ from the second equation of (\ref{homo-equa}). If we write $\hat{g}_n$ into the Fourier series of the type:
\begin{equation}\label{Fourier1}
\hat{g}_n(x,y,t)=\sum\limits_{(k,l)\in\mathbb{Z}^{m+1}}\hat{g}_n^{(k,l)}(y)e^{i(\langle k,\omega\rangle x+lt)},
\end{equation}
where $\hat{g}_n^{(k,l)}(y)$ is the $(k,l)$-Fourier coefficient of $\hat{g}_n(x,y,t)$ with respect to variable $(x,t)$ and
$$\hat{g}_n^{(k,l)}(y)=\lim\limits_{T\rightarrow\infty}\frac{1}{2\pi T}\int_0^T\int_0^{2\pi}\hat{g}_n(x,y,t)e^{-i(\langle k,\omega\rangle x+lt)}dtdx.$$
By (\ref{estimate6}), we have $\hat{g}_n^{(0,0)}(y)=0$.
Assume that
\begin{equation}\label{Fourier2}
v^*(x,y,t)=\sum\limits_{(k,l)\in\mathbb{Z}^{m+1}}v^*_{(k,l)}(y)e^{i(\langle k,\omega\rangle x+lt)},
\end{equation}

Substituting (\ref{Fourier1}) and (\ref{Fourier2}) into the second equation of (\ref{homo-equa}) and comparing the coefficient of $e^{i(\langle k,\omega\rangle x+lt)}$, it follows that for $(k,l)\in\mathbb{Z}^{m+1}\setminus\{(0,0)\}$,
$$v^*_{(k,l)}(y)=i\frac{\hat{g}_n^{(k,l)}(y)}{\langle k,\omega\rangle\gamma+l},$$
where $v^*_{(0,0)}(y)$ will determined later.

Now we write $\hat{f}_n$ into the Fourier series of the type:
\begin{equation}\label{Fourier3}
\hat{f}_n(x,y,t)=\sum\limits_{(k,l)\in\mathbb{Z}^{m+1}}\hat{f}_n^{(k,l)}(y)e^{i(\langle k,\omega\rangle x+lt)}.
\end{equation}
In order to solve the first equation of (\ref{homo-equa}), we have to choose $v^*_{(0,0)}(r)$ such that
\begin{equation}\label{average}
 -v^*_{(0,0)}(y)+\hat{f}_n^{(0,0)}(y)=0,
\end{equation}
so $\hat{f}_n^{(0,0)}(y)=v^*_{(0,0)}(y)$.

 In the same way, we derive
\begin{equation}\label{Fourier4}
u^*(x,y,t)=\sum\limits_{(k,l)\in\mathbb{Z}^{m+1}\setminus\{(0,0)\}}u^*_{(k,l)}(y)e^{i(\langle k,\omega\rangle x+lt)},
\end{equation}
where $u^*_{(0,0)}(y)=0$ and $$u^*_{(k,l)}(y)=i\frac{\hat{f}_n^{(k,l)}(y)-v^*_{(k,l)}(y)}{\langle k,\omega\rangle\gamma+l},\quad (k,l)\in\mathbb{Z}^{m+1}\setminus\{(0,0)\}.$$

$(\mathrm{II})$~ Estimate $u^*(\theta,r,t)$ and $v^*(\theta,r,t)$.

Before we give the estimates of $u^*(\theta,r,t)$ and $v^*(\theta,r,t)$, some technical preparations have to be made.

\begin{lemma}\label{technical1}(Lemma 3.11 in \cite{HR})
Let $\bar{\omega}=(\omega_1,\omega_2,...,\omega_q)\in\mathbb{R}^q$ satisfying the inequalities $D(k,\bar{\omega})\geq\psi(|k|)$, where $D(k,\bar{\omega})=\min\left|\langle(k,l),\bar{\omega}\rangle\right|,k\in\mathbb{Z}^{q-1}\setminus\{0\}$, $\psi$ is an approximation function. Then for $\nu=1,2,...,$ we have
$$\sum\limits_{\substack{\bar{k}\in\mathbb{Z}^q\\0<|\bar{k}|\leq \nu}}\frac{1}{|\langle\bar{k},\bar{\omega}\rangle|^2}\leq\frac{\pi^2}{8}\frac{3^{q+2}}{\psi^2(\nu)},$$
where $\bar{k}=(k,l)$.
\end{lemma}

Choose $$\bar{\omega}=(\omega_1\gamma,\omega_2\gamma,...,\omega_m\gamma,1), q=m+1, \psi(t)=c_0t^{-\sigma} ,$$
 $$k=(k_1,k_2,...,k_m)\in\mathbb{Z}^m\setminus\{0\}, l\in\mathbb{Z},\bar{k}=(k,l).$$
By Lemma \ref{technical1} and the Diophantine conditions (\ref{dioph2}), we obtain
$$\sum\limits_{\substack{(k,l)\in\mathbb{Z}^{m+1}\\0<|k|+|l|\leq \nu}}\frac{1}{|\langle k,\omega\rangle\gamma+l|^2}\leq\frac{\pi^2}{8}3^{m+3}c_0^{-2}\nu^{2\sigma}\leq C\nu^{2\sigma}.$$
\begin{lemma}\label{technical2}
Assume $f(x,t):\{(x,t)\in\mathbb{C}\times\mathbb{T}:|Imx|<s,|Imt|<s\}\rightarrow\mathbb{C}$ be a real analytic and quasi-periodic function with the frequency $\omega=(\omega_1,\omega_2,...,\omega_m)$. Then for $k=(k_1,k_2,...k_m)\in\mathbb{Z}^{m},l\in\mathbb{Z}$, we have the estimate
$$\sum\limits_{(k,l)\in\mathbb{Z}^{m+1}}|f_{(k,l)}|^2e^{2s(|k|+|l|)}\leq 2^{m+1}\|f\|_{s}^2,$$
where
$$f_{(k,l)}=\frac{1}{(2\pi)^{m+1}}\int_{\mathbb{T}^{m}}\int_0^{2\pi}F(\theta,t)e^{-i(\langle k,\theta\rangle+lt)}dtd\theta.$$
are the Fourier coefficients of $f$, $F$ is the shell function of $f(x,\cdot)$.
\end{lemma}
{\bf Proof}
Let $\bar{k}=(k,l)\in\mathbb{Z}^{m+1}$ and $\bar{\theta}=(\theta,t)\in\mathbb{C}^{m+1}$, we have $$\sum\limits_{(k,l)\in\mathbb{Z}^{m+1}}|f_{(k,l)}|^2e^{2s(|k|+|l|)}=\sum\limits_{\bar{k}\in\mathbb{Z}^{m+1}}|f_{\bar{k}}|^2e^{2s|\bar{k}|}$$
and $$f_{(k,l)}=f_{\bar{k}}=\frac{1}{(2\pi)^{m+1}}\int_{\mathbb{T}^{m+1}}F(\bar{\theta})e^{-i\langle \bar{k},\bar{\theta}\rangle}d\bar{\theta}.$$

Set $\bar{x}=(x,t)$ and for every $\lambda=(\lambda_1,\lambda_2,...,\lambda_m,\lambda_{m+1})\in\mathbb{R}^{m+1}$ with $|\lambda|=\max\limits_{1\leq j\leq m+1}|\lambda_j|<s$, define function $\bar{x}+i\lambda\mapsto f(\bar{x}+i\lambda)$,  which domain is $s-|\lambda|$, and its  Fourier coefficient are
$$f_{\bar{k}}(\lambda)=\frac{1}{(2\pi)^{m+1}}\int_{\mathbb{T}^{m+1}}F(\bar{\theta}+i\lambda)e^{-i\langle \bar{k},\bar{\theta}\rangle}d\bar{\theta}.$$
From Bessel's inequality,
$$\sum\limits_{\bar{k}\in\mathbb{Z}^{m+1}}|f_{\bar{k}}(\lambda)|^2\leq\frac{1}{(2\pi)^{m+1}}\int_{\mathbb{T}^{m+1}}|F(\bar{\theta}+i\lambda)|^2 d\bar{\theta}.$$
Hence
\begin{equation}\label{estimate7}
\sum\limits_{\bar{k}\in\mathbb{Z}^{m+1}}|f_{\bar{k}}(\lambda)|^2\leq\|f\|_s^2, \quad|\lambda|<s.
\end{equation}

Define a new function
$$\lambda\mapsto f_{\bar{k}}(\lambda)e^{\langle \bar{k},\lambda\rangle}=\frac{1}{(2\pi)^{m+1}}\int_{\mathbb{T}^{m+1}}F(\bar{\theta}+i\lambda)e^{-i\langle \bar{k},\bar{\theta}+i\lambda\rangle}d\bar{\theta},$$
then
$$\frac{\partial }{\partial \lambda_j}(f_{\bar{k}}(\lambda)e^{\langle \bar{k},\lambda\rangle})=\frac{1}{(2\pi)^{m+1}}\int_{\mathbb{T}^{m+1}}i\frac{\partial }{\partial \bar{\theta}_j}\Big(F(\bar{\theta}+i\lambda)e^{-i\langle \bar{k},\bar{\theta}+i\lambda\rangle}\Big)d\bar{\theta},\quad j=1,2,...,m+1.$$
Since $F(\bar{\theta}+i\lambda)e^{-i\langle \bar{k},\bar{\theta}+i\lambda\rangle}$ is $2\pi$-periodic in $\bar{\theta}_j(j=1,2,...,m+1)$, it follows
$$\frac{\partial }{\partial \lambda_j}(f_{\bar{k}}(\lambda)e^{\langle \bar{k},\lambda\rangle})=0.$$
Thus the function $f_{\bar{k}}(\lambda)e^{\langle \bar{k},\lambda\rangle}$ is independent of $\lambda$, and
$$f_{\bar{k}}=f_{\bar{k}}(0)=f_{\bar{k}}e^{\langle \bar{k},\lambda\rangle}.$$
Consequently, $$|f_{\bar{k}}(\lambda)|^2e^{2\langle \bar{k},\lambda\rangle}=|f_{\bar{k}}|^2,$$
and according to (\ref{estimate7}), we have
\begin{equation}\label{estimate8}
\sum\limits_{\bar{k}\in\mathbb{Z}^{m+1}}|f_{\bar{k}}|^2e^{-2\langle \bar{k},\lambda\rangle}=\sum\limits_{\bar{k}\in\mathbb{Z}^{m+1}}|f_{\bar{k}}(\lambda)|^2\leq\|f\|_s^2, \quad|\lambda|<s.
\end{equation}

Define $e_i\in\mathbb{R}^{m+1}(i=1,2,...,2^{m+1})$ which have components $\pm1$, and
$$\mathbb{Z}_i=\{\bar{k}\in\mathbb{Z}^{m+1}:\langle\bar{k},e_i\rangle=-|\bar{k}|\}.$$
It is easy to see that
\begin{equation}\label{estimate9}
\bigcup\limits_{i=1}^{2^{m+1}}\mathbb{Z}_{i}=\mathbb{Z}^{m+1}.
\end{equation}

Choose $\lambda=\zeta e_i(i=1,2,...,2^{m+1})$ in (\ref{estimate8}), we have
\begin{equation}
\sum\limits_{\bar{k}\in\mathbb{Z}^{i}}|f_{\bar{k}}|^2e^{2\zeta|\bar{k}|}\leq\|f\|_s^2, \quad 0<\zeta<s.
\end{equation}
Passing to the limit $\zeta\rightarrow s$ yields
\begin{equation}\label{estimate10}
\sum\limits_{\bar{k}\in\mathbb{Z}^{i}}|f_{\bar{k}}|^2e^{2s|\bar{k}|}\leq\|f\|_s^2.
\end{equation}
Combing (\ref{estimate9}) with (\ref{estimate10}), we have
\begin{equation}\label{estimate11}
\sum\limits_{\bar{k}\in\mathbb{Z}^{m+1}}|f_{\bar{k}}|^2e^{2s|\bar{k}|}\leq 2^{m+1}\|f\|_s^2.
\end{equation}
We complete the proof of this lemma.$\hfill\square$

Since function $\hat{g}_n(x,y,t)$ is real analytic and quasi-periodic with the frequency $\omega=(\omega_1,\omega_2,...,\omega_m)$, as an application of Lemma \ref{technical2} to $\hat{g}_n(x,y,t)$ on $D(s_n,r_n)$ yields that
\begin{equation}\label{estimate12}
\sum\limits_{(k,l)\in\mathbb{Z}^{m+1}\setminus\{(0,0)\}}|\hat{g}_n^{(k,l)}(y)|^2e^{2s_n(|k|+|l|)}\leq C\|\hat{g}_n\|_{s_n,r_n}^2.
\end{equation}

In the following, we will estimate  $u,v,u^*$ and $v^*$.

Firstly, we estimate the sum
$$G_\nu(y)=\sum\limits_{\substack{(k,l)\in\mathbb{Z}^{m+1}\\0<|k|+|l|\leq \nu}}\left|\frac{\hat{g}_n^{(k,l)}(y)}{\langle k,\omega\rangle\gamma+l}\right|e^{(|k|+|l|)s_n},\quad \nu=1,2,\cdots.$$

In view of Cauchy-Schwarz inequality and Lemma \ref{technical1}, we have
\begin{equation*}
\begin{aligned}
G_\nu(y)\leq&\sqrt{\sum\limits_{\substack{(k,l)\in\mathbb{Z}^{m+1}\\0<|k|+|l|\leq \nu}}\left|\hat{g}_n^{(k,l)}(y)\right|^2e^{2(|k|+|l|)s_n}}\cdot\sqrt{\sum\limits_{\substack{(k,l)\in\mathbb{Z}^{m+1}\\0<|k|+|l|\leq \nu}}|\langle k,\omega\rangle\gamma+l|^{-2}}\\
\leq&C\|\hat{g}_n\|_{s_n,r_n}\nu^\sigma.
\end{aligned}
\end{equation*}
Set $G_0(y)=0$, we obtain
\begin{equation*}
\sum\limits_{\substack{(k,l)\in\mathbb{Z}^{m+1}\\0<|k|+|l|\leq N}}\left|\frac{\hat{g}_n^{(k,l)}(y)}{\langle k,\omega\rangle\gamma+l}\right|e^{(|k|+|l|)(s_n-\rho_n)}=(1-e^{-\rho_n})\sum\limits_{\nu=1}^{N}G_\nu(y)e^{-\nu\rho_n}+G_N(y)e^{-(N+1)\rho_n},
\end{equation*}
where $\rho_n=\frac{1}{200p}(s_n-s_{n+1})$.

Taking $N\rightarrow\infty$, we derive
\begin{equation*}
\begin{aligned}
&\sum\limits_{(k,l)\in\mathbb{Z}^{m+1}\setminus\{(0,0)\}}\left|\frac{\hat{g}_n^{(k,l)}(y)}{\langle k,\omega\rangle\gamma+l}\right|e^{(|k|+|l|)(s_n-\rho_n)}\\
&\leq(1-e^{-\rho_n})\sum\limits_{\nu=1}^{\infty}G_\nu(y)e^{-\nu\rho_n}\\
&\leq C\|\hat{g}_n\|_{s_n,r_n}\sum\limits_{\nu=1}^{\infty} \nu^\sigma\Big(e^{-\nu\rho_n}-e^{-(\nu+1)\rho_n}\Big)\\
&\leq C\|\hat{g}_n\|_{s_n,r_n}\rho_n^{-\sigma}\leq C\|\hat{g}_n\|_{s_n,r_n}s_n^{-\sigma}.
\end{aligned}
\end{equation*}
Hence, combing with (\ref{estimate5}),(\ref{Fourier2}) and (\ref{average}) , it follows that
\begin{equation}\label{estimate13}
\begin{aligned}
&\|v^*(x,y,t)\|_{s_n^{1},r_n}\leq C\|\hat{g}_n\|_{s_n,r_n}s_n^{-\sigma}\\
&\leq C\varepsilon\cdot\varepsilon_{n-1}s_{n-1}^{m}s_n^{-\sigma}\leq C\varepsilon\cdot\tau_{n}\varepsilon_{n}s_{n}^{-\frac{\mu}{100}}.
\end{aligned}
\end{equation}

According to the same process, we obtain
\begin{equation*}
\sum\limits_{(k,l)\in\mathbb{Z}^{m+1}\setminus\{(0,0)\}}|\hat{f}_n^{(k,l)}(y)|^2e^{2s_n(|k|+|l|)}\leq C\|\hat{f}_n\|_{s_n,r_n}^2
\end{equation*}
and
\begin{equation}\label{estimate14}
\begin{aligned}
&\|u^*(x,y,t)\|_{s_n^{2},r_n}\leq C\|v^*(x,y,t)\|_{s_n,r_n}s_n^{-\sigma}+C\|\hat{f}_n(x,y,t)\|_{s_n,r_n}s_n^{-\sigma}\\
&\leq C\varepsilon\cdot\tau_{n}\varepsilon_{n}s_{n}^{-(\sigma+\frac{\mu}{100})}.
\end{aligned}
\end{equation}
Using Cauchy estimate on the derivatives of $u^*(x,y,t)$ and $v^*(x,y,t)$, we get the following estimates:
\begin{equation}\label{estimates15}
\|\partial_x^k\partial_y^lu^*\|_{s_n^2,r_n^1}\leq C\varepsilon\cdot\tau_{n}\varepsilon_{n}s_{n}^{-(\sigma+\frac{\mu}{100})}\max\{r_n^{-1},s_n^{-1}\},
\end{equation}
\begin{equation}\label{estimates16}
\|\partial_x^k\partial_y^lv^*\|_{s_n^2,r_n^1}\leq C\varepsilon\cdot\tau_{n}\varepsilon_{n}s_{n}^{-\frac{\mu}{100}}\max\{r_n^{-1},s_n^{-1}\},
\end{equation}
where $k+l=1$, $k\geq0$ and $l\geq 0$.

By the implicit function theorem, we know that $\Delta\Phi_{n+1}=\Phi=\Psi^{-1}:$
\begin{equation}
\begin{cases}
x=\xi+u(\xi,\eta,t),\\
y=\eta+v(\xi,\eta,t),
\end{cases}
\end{equation}
where $u,v$ are real analytic, and by Lemma \ref{quasi-prop}, it implies that $u,v$ are quasi-periodic with the frequency $\omega=(\omega_1,\omega_2,...,\omega_m)$. Moreover,
\begin{align}
&\|u(\xi,\eta,t)\|_{s_n^{4},r_n^2}\leq C\varepsilon\cdot\tau_{n}\varepsilon_{n}s_{n}^{-(m+\frac{11\mu}{100})},\label{estimates17}\\
&\|v(\xi,\eta,t)\|_{s_n^{4},r_n^2}\leq C\varepsilon\cdot\tau_{n}\varepsilon_{n}s_{n}^{-\frac{\mu}{100}},\label{estimates17-1}
\end{align}
and
$$\Delta\Phi_{n+1}\big(D(s_{n+1},r_{n+1})\big)\subset\Delta\Phi_{n+1}\big(D(s_{n}^4,r_{n}^2)\big)\subset D(s_{n},r_{n}).$$

Next, we are in position to prove the change $\Delta\Phi_{n+1}$ is $G-$invariant with respect to $G(\xi,\eta,t)=(-\xi,\eta,-t)$.

By means of (\ref{estimate6}), (\ref{Fourier2}) and (\ref{Fourier4}), we have
$$u^*(-x,y,-t)=-u^*(x,y,t),\quad v^*(-x,y,-t)=v^*(x,y,t).$$
Combing with $\Phi\circ\Psi=id$, it follows that
$$u^*(x,y,t)+u\big(x+u^*(x,y,t),y+v^*(x,y,t),t\big)=0,$$
and $$v^*(x,y,t)+v\big(x+u^*(x,y,t),y+v^*(x,y,t),t\big)=0.$$

By the discussion as above, it yields that
\begin{equation}\label{even-odd1}
u(-\xi,\eta,-t)=-u(\xi,\eta,t),\quad v(-\xi,\eta,-t)=v(\xi,\eta,t),
\end{equation}
where $(\xi,\eta,t)\in D(s_{n}^4,r_{m}^2)$. Therefore the transformation $\Delta\Phi_{n+1}$ is $G-$invariant with respect to $G(\xi,\eta,t)=(-\xi,\eta,-t)$.

$(\mathrm{III})$~ Estimate the new  nonlinear parts.

From (\ref{change4}) and (\ref{change5}), the new perturbations $\bar{f}_{n+1}(\xi,\eta,t)$ and $\bar{g}_{n+1}(\xi,\eta,t)$ are
\begin{align}
&\partial _y u^*\cdot\hat{g}_n(x,y,t)+\partial _x u^*\cdot\hat{f}_n(x,y,t)\label{new-pertu1-1}\\
&+\partial _x u^*\cdot y,\label{new-pertu1-2}\\
&\partial _y v^*\cdot\hat{g}_n(x,y,t)+\partial _x v^*\cdot\hat{f}_n(x,y,t)\label{new-pertu2-1}\\
&+\partial _xv^*\cdot y,\label{new-pertu2-2}
\end{align}
repectively. observe $\bar{f}_{n+1}(\xi,\eta,t)$ and $\bar{g}_{n+1}(\xi,\eta,t)$ and by means of Lemma \ref{quasi-prop}, it is obvious that  $\bar{f}_{n+1}(\xi,\eta,t)$ and $\bar{g}_{n+1}(\xi,\eta,t)$ are real analytic and quasi-periodic with the frequency $\omega=(\omega_1,\omega_2,...,\omega_m)$.

First of all, we estimate $\bar{f}_{n+1}(\xi,\eta,t)$.

For (\ref{new-pertu1-1}), by (\ref{estimate5}), (\ref{estimates15})  and regarding (\ref{new-pertu1-1}) as a function of $(x,y,t)$, we have
$$\|(\ref{new-pertu1-1})(x,y,t)\|_{s_{n}^4,r_{n}^2}\leq \frac{1}{2}C\varepsilon\cdot\varepsilon_{n}.$$

For (\ref{new-pertu1-2}), by Cauchy estimate and in view of (\ref{estimates15}), we have
\begin{equation*}
\begin{aligned}
\|(\ref{new-pertu1-2})(x,y,t)\|_{s_{n+1},r_{n+1}}&\leq \|\partial _x u^*\cdot (\eta-v^*)\|_{s_{n+1},r_{n+1}}\\
&\leq\|\partial _x u^*\cdot\eta\|_{s_{n+1},r_{n+1}}+\|\partial _x u^*\|_{s_{n}^4,r_{n}^2}\|v^*\|_{s_{n}^4,r_{n}^2}\\
&\leq C\varepsilon\cdot\tau_{n}\varepsilon_{n}s_{n}^{-(\sigma+\frac{\mu}{100})}s_n^{-1}r_{n+1}\\
&+C\varepsilon^2\cdot\tau_{n}\varepsilon_{n}s_{n}^{-(\sigma+\frac{\mu}{100})}s_n^{-1}\tau_{n}\varepsilon_{n}s_{n}^{-\frac{\mu}{100}}\\
&\leq\frac{1}{2}C\varepsilon\cdot\varepsilon_{n}.
\end{aligned}
\end{equation*}
Hence,  $$\|\bar{f}_{n+1}(\xi,\eta,t)\|_{s_{n+1},r_{n+1}}\leq C\varepsilon\cdot\varepsilon_{n}.$$

In a similar way, we see that $$\|\bar{g}_{n+1}(\xi,\eta,t)\|_{s_{n+1},r_{n+1}}\leq C\varepsilon\cdot\varepsilon_{n}s_n^{m}.$$

Therefore, we complete the proof of inequality (\ref{estimate2}) with replacing $n$ by $n+1$.

By the definition of  $\bar{f}_{n+1}(\xi,\eta,t), \bar{g}_{n+1}(\xi,\eta,t)$ and combining with (\ref{estimate6}) and (\ref{even-odd1}), we verify that $\bar{f}_{n+1}(\xi,\eta,t)$ and $\bar{g}_{n+1}(\xi,\eta,t)$ meet (\ref{estimate3}) with replacing $n$ by $n+1$.

Consequently, we prove the Iteration Lemma.

\section{The proof of Theorem \ref{theo-map} and Theorem \ref{theo-sys}}

Before we verify the Theorem \ref{theo-map} and Theorem \ref{theo-sys}, we show that there is a convergent change of variables, transforming the given system (\ref{main-sys}) into the linearized normal form
\begin{equation}
\xi'=\gamma+\eta,\quad\eta'=0.
\end{equation}

On the one hand, by the discussion of the previous section, we know there exists a transformation sequence $\{\Phi_n\}$, the convergence of which is decided by their nonlinear parts
$$\bar{u}_n=u_1+u_2+\cdots+u_n,$$
$$\bar{v}_n=v_1+v_2+\cdots+v_n.$$
From (\ref{change-esti1}) and (\ref{change-esti2}), we have $\|\bar{u}_n\|+\|\bar{v}_n\|\rightarrow 0$ as $n\rightarrow0$, it implies that  there exists a subsequence of $\{\Phi_n\}$, which converges to a transformation $\Phi_\infty$. Without loss of generality, one can assume the subsequence is $\{\Phi_n\}$, it follows that$$\lim_{n\rightarrow\infty}\Phi_n=\Phi_\infty.$$

On the other hand, since that $\|\bar{f}_n\|_{s_n,r_n}+ \|\bar{g}_n\|_{s_n,r_n}\rightarrow 0$ as $n\rightarrow0$,
 we have that
\begin{equation}
\begin{aligned}
&(\partial\Phi_\infty)^{-1}\cdot \big(\gamma+y+f(x,y,t),g(x,y,t)\big)\circ\Phi_\infty\\
=&\lim_{n\rightarrow\infty} \big(\partial \Phi_n)^{-1}\cdot (\gamma+y+f(x,y,t),g(x,y,t)\big)\circ\Phi_n\\
=&\lim_{n\rightarrow\infty}\big(\gamma+y+\bar{f}_{n}(\xi,\eta,t),\bar{g}_{n}(\xi,\eta,t)\big)\\
=&(\gamma+y,0).
\end{aligned}
\end{equation}

Let $$\psi(x,t)=\Phi_\infty(x,t)=\lim_{n\rightarrow\infty}\Phi_n(x,0,t),\quad (x,t)\in\mathbb{R}\times\mathbb{T}.$$
It suffices to prove the fact that $\psi(x,t)$ does exist and is $\mathcal{C}^1$. Then $\psi(\mathbb{R}\times\{y=0\}\times\mathbb{T})$ is an invariant tours for the original system (\ref{main-sys}).
Denote
$$\Phi^j=\Phi_j(x,0,t)$$
and write
$$\Phi^n=\Phi^0+\sum_{j=1}^{j=n}(\Phi^j-\Phi^{j-1}).$$

Then by (\ref{change-esti1}) and (\ref{change-esti2}), we have
$$|\psi(x,0,t)-Id(x,0)|\leq|\Phi^0(x,0,t)-Id(x,0)|+\sum_{j=1}^\infty|\Phi^j-\Phi^{j-1}|\leq C\varepsilon\cdot\varepsilon_{0}s_{0}^{-\frac{\mu}{100}}\varepsilon^{-(1+\tilde{\mu})^{-1}\cdot(1+\frac{m}{p})\tilde{\mu}}.$$
This means that $\psi$ does exist and is $\mathcal{C}^0$. Moreover, $\psi$ is quasi-periodic in $x$ with the frequency $\omega=(\omega_1,\omega_2,...,\omega_m)$ and
\begin{equation}\label{conver}
\psi=Id+O(\varepsilon).
\end{equation}

In the following, we are ready to prove $\psi$ is injective by the wonderful idea due to R\"{u}ssmann \cite{HR1}. Let $\tilde{\gamma}=(\gamma,1)$. Although the system (\ref{main-sys}) is non-autonomous, it can be rewritten as an autonomous one:
\begin{equation*}
\begin{cases}
x'=\gamma+y+f(x,y,\theta),\\
y'=g(x,y,\theta),\\
\theta'=1
\end{cases}
(x,y,\theta)\in D.
\end{equation*}
Observe that for any $w=(x,t)\in\mathbb{R}\times\mathbb{T}$, $\varsigma(w,t):=\psi(w+\tilde{\gamma}t)$ is the solution to autonomous differential equation with initial $\psi(w)$. Take $w_1,w_2\in\mathbb{R}\times\mathbb{T}$ and let $\eta=w_1-w_2$. Suppose that $\psi(w_1)=\psi(w_2)$. By Picard's existence and uniqueness theorem,
$$\psi(w_1+\tilde{\gamma}t)=\psi(w_2+\tilde{\gamma}t), \quad t\in\mathbb{R}.$$
Since $\tilde{\gamma}$ is in Diophantine class, the set $\{\tilde{\gamma}t:t\in\mathbb{R}\}$ is dense in $\mathbb{R}\times\mathbb{T}$. So
$$\psi(w_1+a)=\psi(w_2+a), \quad a\in\mathbb{R}\times\mathbb{T}.$$
Choose $a=-w_2$, we have $\psi(0)=\psi(\eta)$. Repeating the procedure as above, we have
$$\psi(l\eta)=\psi(0),\quad l\in\mathbb{Z}^+.$$
From (\ref{conver}), it yields that $\psi(l\eta)=O(\varepsilon)$. Taking $l\rightarrow+\infty$, one has $\eta=0$. Thus $\psi$ is injective. Moreover, $\psi: \mathbb{R}\times\mathbb{T}\rightarrow\psi(\mathbb{R}\times\mathbb{T})\subset\mathbb{R}^2$ is a homomorphism. Therefore, $\psi(\mathbb{R}\times\mathbb{T})$ is an invariant torus with rational frequency vector $\tilde{\gamma}$.

{\bf Proof of theorem \ref{theo-sys}} To sum up, we see that $\psi\big(\mathbb{R}\times \{y=0\}\times\mathbb{T}\big)$ is an invariant torus of the system (\ref{main-sys}).

{\bf Proof of theorem \ref{theo-map}} ~According to Proposition 4.5 of Sevryuk \cite{MBS}, the mapping (\ref{some-resu}) in Theorem \ref{theo-map} can be regarded as the corresponding Poincar\'{e} mapping of system (\ref{main-sys}). The proof is completed by Theorem \ref{theo-sys}.

\section{The proof of Theorem \ref{theo-map2}}

The proof is similar to the \cite[Theorem 1]{LB}. Thus we just give the sketch here.

Since the mapping $M_\delta$ is reversible with respect to the involution $G:(x,y)\rightarrow(-x,y)$ for any $\delta\in[0,1]$, it follows that
\begin{equation}
\begin{cases}
l_1(x,y)+f(x,y;\delta)=l_1\circ G\circ M_\delta(x,y)+f\big(G\circ M_\delta(x,y);\delta\big),\\
$$l_2(x,y)+g(x,y;\delta)=-l_2\circ G\circ M_\delta(x,y)-g\big(G\circ M_\delta(x,y);\delta\big),
\end{cases}
\end{equation}
for any $\delta\in[0,1]$ and $(x,y)\in\mathbb{R}\times B(r_0)$.
Taking $\delta\rightarrow 0$, it yields that
\begin{equation}\label{small-symm}
l_1(x,y)=l_1(-x-\gamma,y),\quad l_2(x,y)=-l_2(-x-\gamma,y),\quad (x,y)\in\mathbb{R}\times B(r_0).
\end{equation}
Writting $l_1$ and $l_2$ into the Fourier series of the type:
$$l_1(x,y)=\sum\limits_{k\in\mathbb{Z}^m}l_1(k,y)e^{i\langle k,\omega\rangle x},\quad l_2(x,y)=\sum\limits_{k\in\mathbb{Z}^m}l_2(k,y)e^{i\langle k,\omega\rangle x},$$
with $l_1(-k,y)=\bar{l}_1(k,y),l_2(-k,y)=\bar{l}_2(k,y)$.

Combing with (\ref{small-symm}), we obtain
\begin{equation}\label{small-symm1}
l_1(k,y)=l_1(-k,y)e^{i\langle k,\omega\rangle\gamma}, \quad l_2(k,y)=-l_2(-k,y)e^{i\langle k,\omega\rangle\gamma}.
\end{equation}
In particular, from the above discussion and our assumption (\ref{average1}), one has $$h(y)=l_1(0,y)=\lim\limits_{T\rightarrow\infty}\frac{1}{T}\int_0^T\frac{\partial l_1}{\partial x}(x,y)dx>0,\quad l_2(0,y)=\lim\limits_{T\rightarrow\infty}\frac{1}{T}\int_0^Tl_2(x,y)dx=0.$$

Let
$$h_1(x,y)=\sum\limits_{\substack{k\in\mathbb{Z}^m\\0<|k|< N}}l_1(k,y)e^{i\langle k,\omega\rangle x},\quad h_2(x,y)=\sum\limits_{\substack{k\in\mathbb{Z}^m\\0<|k|<N}}l_2(k,y)e^{i\langle k,\omega\rangle x}.$$
It is well known that for any $\varepsilon'>0$, there exists a positive integer $N$ depending on $\varepsilon', l_1, l_2$ such that
\begin{equation}\label{small-symm2}
\|l_1(x,y)-h(y)-h_1(x,y)\|_{\mathcal{C}^{p+m}}+\|l_2(x,y)-h_2(x,y)\|_{\mathcal{C}^{p+m}}<\varepsilon'.
\end{equation}

In what follows, we consider the difference equations:
\begin{equation}
\begin{cases}
u(x+\gamma,y)-u(x,y)+h_1(x,y)=0,\\
v(x+\gamma,y)-v(x,y)+h_2(x,y)=0.
\end{cases}
\end{equation}
It is easy to derive that
$$u(x,y)=-\sum\limits_{\substack{k\in\mathbb{Z}^m\\0<|k|< N}}\frac{l_1(k,y)}{e^{i\langle k,\omega\rangle\gamma}-1}e^{i\langle k,\omega\rangle x},$$
$$v(x,y)=-\sum\limits_{\substack{k\in\mathbb{Z}^m\\0<|k|< N}}\frac{l_2(k,y)}{e^{i\langle k,\omega\rangle\gamma}-1}e^{i\langle k,\omega\rangle x}.$$
Moreover, by (\ref{small-symm1}), we have
$$u(-x,y)=-u(x,y),\quad v(-x,y)=v(x,y).$$
Let $\zeta_0(N)=2N\max\limits_{0<k\leq N}|e^{i\langle k,\omega\rangle\gamma}-1|^{-1}$, then
$$\|u\|_{\mathcal{C}^{p+m}}+\|v\|_{\mathcal{C}^{p+m}}\leq \zeta_0(N)(\|l_1\|_{\mathcal{C}^{p+m}}+\|l_2\|_{\mathcal{C}^{p+m}}).$$

Let $$R_1(x,y)=u(x+\gamma,y)-u(x,y)+l_1(x,y)-h(y),\quad R_2(x,y)=v(x+\gamma,y)-v(x,y)+l_2(x,y),$$
then by (\ref{small-symm2}), one has that
$$\|R_1\|_{\mathcal{C}^{p+m}}+\|R_2\|_{\mathcal{C}^{p+m}}<\varepsilon'.$$

Define the change of variables $\Phi$ by
\begin{equation}\label{small-symm3}
\begin{cases}
\xi=x+\delta u(x,y),\\
\eta=y+\delta v(x,y).
\end{cases}
\end{equation}
Then the mapping $M_\delta$ can be transformed into $L_\delta=\Phi\circ M_\delta\circ\Phi^{-1}:$
\begin{equation}
\begin{cases}
\xi_1=\xi+\gamma+\delta h(\eta)+\delta \varphi_1\circ\Phi^{-1} (\xi,\eta;\delta),\\
\eta_1=\eta+\delta \varphi_2 \circ\Phi^{-1}(\xi,\eta;\delta),
\end{cases}
\end{equation}
where
$$\varphi_1(x,y;\delta)=R_1(x,y)+u(x_1,y_1)-u(x+\gamma,y)+h(y)-h(\eta)+f(x,y;\delta),$$
$$\varphi_2(x,y;\delta)=R_2(x,y)+v(x_1,y_1)-v(x+\gamma,y)+g(x,y;\delta).$$

Note that $l_1,l_2,f,g$ are quasi-periodic in $x$ with the frequency $\omega=(\omega_1,\omega_2,...,\omega_m)$, by the definition of $\varphi_1,\varphi_2$,  $\varphi_1,\varphi_2$ are also quasi-periodic in $x$ with the frequency $\omega=(\omega_1,\omega_2,...,\omega_m)$. Hence $\varphi_1\circ\Phi^{-1},\varphi_1\circ\Phi^{-1}$ are also quasi-periodic in $x$ with the frequency $\omega$, which is guaranteed by the definition $\Phi$ and Lemma \ref{quasi-prop}. Since $M_\delta$ is reversible with respect to $G(x,y)\rightarrow(-x,y)$ and $\Phi$ is $G$-invariant, it means that $L_\delta$ is reversible with respect to $G(\xi,\eta)\rightarrow(-\xi,\eta)$.

In the following, we are in position to estimate $\varphi_1\circ\Phi^{-1},\varphi_1\circ\Phi^{-1}$.
Similar to \cite{LB}, there exists a constant $C_1>0$ such that
$$\|\varphi_1\circ\Phi^{-1}\|_{\mathcal{C}^p}+\|\varphi_2\circ\Phi^{-1}\|_{\mathcal{C}^{p+m}}\leq C_1(\|\varphi_1\|_{\mathcal{C}^p}+\|\varphi_2\|_{\mathcal{C}^{p+m}}).$$

From Lemma 2.2 in \cite{LB}, it follows that

\begin{equation}
\begin{aligned}
&\|u(x_1,y_1)-u(x+\gamma,y)\|_{\mathcal{C}^{p}}+\|v(x_1,y_1)-v(x+\gamma,y)\|_{\mathcal{C}^{p+m}}\\
&\leq C_2\delta \big(\|l_1\|_{\mathcal{C}^{p+m}}+\|l_2\|_{\mathcal{C}^{p+m}}+\|f\|_{\mathcal{C}^{p}}+\|g\|_{\mathcal{C}^{p+m}}\big)\\
&\leq C_2\delta\big(\|l_1\|_{\mathcal{C}^{p+m}}+\|l_2\|_{\mathcal{C}^{p+m}}+1\big).
\end{aligned}
\end{equation}

Since $f(x,y;0)=g(x,y;0)=0$, there is a $\Delta_1>0$ such that for $\delta\in(0,\Delta_1)$,
$$\|f\|_{\mathcal{C}^{p}}+\|g\|_{\mathcal{C}^{p+m}}<\varepsilon.$$

According to the above discussion, it follows that
$$\|\varphi_1\|_{\mathcal{C}^p}+\|\varphi_2\|_{\mathcal{C}^{p+m}}\leq C_1\Big(\varepsilon+\varepsilon'+C_2\delta \big(\|l_1\|_{\mathcal{C}^{p+m}}+\|l_2\|_{\mathcal{C}^{p+m}}+1\big)+\Omega(\delta)\Big).$$

Choose $$\varepsilon<\frac{1}{4C_1}\varepsilon_0,\quad\varepsilon'<\frac{1}{4C_1}\varepsilon_0,$$ where $\varepsilon_0$ is given in (\ref{perturbation1}).

In view of (\ref{average1}), we derive $h'(y)>0$. By (\ref{small-symm3}), one has that $\|h(\eta)-h(y)\|=\|h(y+\delta v)-h(y)\|\leq \Omega(\delta)$, where $\Omega(\delta)$ is an appropriate modulus of continuity that depends on $l_1$ and $v$. Due to that $\lim\limits_{\delta\rightarrow 0}\Omega(\delta)=0$, it means that there is a constant $\Delta_2$ such that for $\delta\in(0,\Delta_2)$, $$\Omega(\delta)<\frac{1}{4C_1}\varepsilon_0.$$

Let $$\bar{\Delta}=\min\{\Delta_0,\Delta_1,\Delta_2,\frac{1}{4C_1C_2\big(\|l_1\|_{\mathcal{C}^{p+m}}+\|l_2\|_{\mathcal{C}^{p+m}}+1\big)}\varepsilon_0\},$$
where $\Delta_0$ is given in Remark \ref{zhu},  Then for $\delta\in(0,\bar{\Delta})$, we have
\begin{equation}
\begin{aligned}
&\|\varphi_1\circ\Phi^{-1}\|_{\mathcal{C}^p}+\|\varphi_2\circ\Phi^{-1}\|_{\mathcal{C}^{p+1}}\\
&\leq C_1(\|\varphi_1\|_{\mathcal{C}^p}+\|\varphi_2\|_{\mathcal{C}^{p+m}})\\
&\leq\varepsilon_0.
\end{aligned}
\end{equation}
Thus, for $\delta\in(0,\bar{\Delta})$, $$\|\varphi_1\circ\Phi^{-1}\|_{\mathcal{C}^p}\leq\varepsilon_0,\quad\|\varphi_2\circ\Phi^{-1}\|_{\mathcal{C}^{p+m}}\leq\varepsilon_0.$$

Therefore the mapping $L_\delta$ meets all the assumptions of Remark \ref{zhu}, and $L_\delta$ has invariant curves. Undoing the change of variables, we derive that the existence of invariant curves of mapping $M_\delta$. We complete the proof of Theorem \ref{theo-map2}.

\section{Application}
In this section, we will apply Theorem \ref{theo-map} to the equation
 \begin{equation}\label{appl-equa}
x''+\varphi(x)f(x')+\omega^2x+g(x)=p(t).
 \end{equation}

Suppose that

$(H_1)$: $f,g,\varphi\in\mathcal{C}^{3m+3}$ and $p\in\mathcal{C}^{4m+5}$ ;

$(H_2)$: $f$ and $p$ are even functions, $p(t)$ are quasi-periodic in $t$ with the frequency $\mu=(\mu_1,\mu_2,...,\mu_m)$;

$(H_3)$: $\lim_{x\rightarrow\pm\infty}\varphi(x)=:\varphi(\pm\infty)\in \mathbb{R},\quad \lim_{|x|\rightarrow+\infty}x^{3m+3}\varphi^{(3m+3)}(x)=0;$

$(H_4)$: $\lim_{x\rightarrow\pm\infty}f(x)=:f(+\infty)\in \mathbb{R},\quad \lim_{|x|\rightarrow+\infty}x^{3m+3}f^{(3m+3)}(x)=0;$

$(H_5)$: $\lim_{x\rightarrow\pm\infty}g(x)=:g(\pm\infty)\in \mathbb{R},\quad \lim_{|x|\rightarrow+\infty}x^{3m+3}g^{(3m+3)}(x)=0.$

\begin{theorem}\label{appl-theo}
Suppose that $(H_1)-(H_5)$ hold, and
$\omega$ satisfies
\begin{equation}\label{diophantine}
|\langle k,\mu\rangle\omega^{-1}-l|\geq\frac{c_0}{|k|^{\sigma}},
\end{equation}
where $c_0>0,\sigma>0, k\in\mathbb{Z}\setminus 0, l\in\mathbb{Z}$.
Then for every solution $x(t)$ of (\ref{appl-equa}), we have
$$\sup_{t\in\mathbb{R}}(|x(t)|+|x'(t)|)<+\infty.$$
\end{theorem}

In order to obtain the boundedness of all solutions of (\ref{appl-equa}), it is sufficient to prove that its Poincar\'{e} mapping can be written as a twist mapping with small enough perturbations. Under some transformations, we apply Theorem \ref{theo-map} to achieve the goal.  In the following, we will give the proof of Theorem \ref{appl-theo}, which is similar to the proof in \cite{LB3} and \cite{Li2}. Thus we give the sketch of the proof.

We first rewrite (\ref{appl-equa}) as
\begin{equation}\label{appl-syst1}
\begin{cases}
x'=-\omega y\\
y'=\omega x+\omega^{-1}\varphi(x)f(\omega y)+\omega^{-1}g(x)-\omega^{-1}p(t).
\end{cases}
\end{equation}
From $(H_2)$, it follows that (\ref{appl-syst1}) is reversible with respect to the involution $G(x,y)=(x,-y)$.

By polar coordinates change
$x=r\cos\theta,\ y=r\sin\theta,$
the system (\ref{appl-syst2}) is transformed into
\begin{equation}\label{appl-syst2}
\begin{cases}
r'=\omega^{-1}\big(\varphi(r\cos\theta)f(\omega r\sin\theta)+g(r\cos\theta)\big)\sin\theta-\omega^{-1}p(t)\sin\theta\\
\theta'=\omega +\omega^{-1}r^{-1}\big(\varphi(r\cos\theta)f(\omega r\sin\theta)+g(r\cos\theta)\big)\cos\theta-\omega^{-1}r^{-1}p(t)\cos\theta.
\end{cases}
\end{equation}
Observing that
$$|\omega^{-1}r^{-1}\big(\varphi(r\cos\theta)f(\omega r\sin\theta)+g(r\cos\theta)\big)\cos\theta-\omega^{-1}r^{-1}p(t)\cos\theta|\leq Cr^{-1},$$
for some $C>0$, we may consider (\ref{appl-syst2}) assuming that $r(t)>2C\omega^{-1}$ for all $t\in\mathbb{R}$ along a solution $t\mapsto(r(t),\theta(t))$. Therefore,
$$\theta'\geq\frac{1}{2}\omega>0,\quad t\in\mathbb{R},$$
which means that $t\mapsto\theta(t)$ ia globally invertible. Denoting by $\theta\mapsto t(\theta)$ the inverse function, we have that $\theta\mapsto (r(t(\theta)),t(\theta))$ solves the system
\begin{equation}\label{appl-syst3}
\begin{cases}
\frac{dr}{d\theta}=\Phi(r,t,\theta)\\
\frac{dt}{d\theta}=\Psi(r,t,\theta),
\end{cases}
\end{equation}
where $$\Phi(r,t,\theta)=\frac{\omega^{-1}\big(\varphi(r\cos\theta)f(\omega r\sin\theta)+g(r\cos\theta)\big)\sin\theta-\omega^{-1}p(t)\sin\theta}{\omega +\omega^{-1}r^{-1}\big(\varphi(r\cos\theta)f(\omega r\sin\theta)+g(r\cos\theta)\big)\cos\theta-\omega^{-1}r^{-1}p(t)\cos\theta},$$
$$\Psi(r,t,\theta)=\frac{1}{\omega +\omega^{-1}r^{-1}\big(\varphi(r\cos\theta)f(\omega r\sin\theta)+g(r\cos\theta)\big)\cos\theta-\omega^{-1}r^{-1}p(t)\cos\theta}.$$
Now noting that the action, angle and time variables are $r, t$ and $\theta$, respectively. Since $\Psi(r,-t,-\theta)=\Psi(r,t,\theta)$ and $\Phi(r,-t,-\theta)=-\Phi(r,t,\theta)$, we see that system (\ref{appl-syst3}) is reversible under the transformation $(r,t)\mapsto(r,-t)$.

To estimate error terms, we introduce some notations.
\begin{definition}\label{defin} \item[{\textrm{(i)}}]: Assume function $f(\theta,r,t)$ is $O_n(r^{-j})$, if $f$ is smooth in $(r,t)$, continue in $\theta$, periodic of period $2\pi$ in $\theta$ and $t$, moreover
$$|r^{k+j}\frac{\partial^{k+l}f}{\partial r^k\partial t^l}|\leq C,\quad 0\leq k+l\leq n,$$
where $C$ is a positive constant.

\item[{\textrm{(ii)}}]: Suppose function $f(\theta,r,t)$ is $o_n(r^{-j})$, if $f$ is smooth in $(r,t)$, continue in $\theta$, periodic of period $2\pi$ in $\theta$ and $t$, moreover
$$\lim_{r\rightarrow\infty}r^{k+j}\frac{\partial^{k+l}f}{\partial r^k\partial t^l}=0,\quad 0\leq k+l\leq n,$$
uniformly in $(\theta,t)$.
\end{definition}
It is obvious that $$\Phi(r,t,\theta)\in O_{3m+3}(1),\quad\Psi(r,t,\theta)\in O_{3m+3}(1), $$
and (\ref{appl-syst3}) can be rewritten as
\begin{equation}\label{appl-syst4}
\begin{cases}
\frac{dr}{d\theta}=\omega^{-2}\big(\varphi(r\cos\theta)f(\omega r\sin\theta)+g(r\cos\theta)\big)\sin\theta-\omega^{-2}p(t)\sin\theta+O_{3m+3}(r^{-1})\\
\frac{dt}{d\theta}=\omega^{-1}-\omega^{-3}r^{-1}\big(\varphi(r\cos\theta)f(\omega r\sin\theta)+g(r\cos\theta)\big)\cos\theta+\omega^{-3}r^{-1}p(t)\cos\theta+O_{3m+3}(r^{-2}).
\end{cases}
\end{equation}

Since the Poincar\'{e} mapping of (\ref{appl-syst4}) is not sufficiently close to a twist map, we need to transform (\ref{appl-syst4}) further.

Let
$$\lambda=r+S_1(\theta,r),\quad t=t,$$
where
$$S_1(\theta,r)=- \omega^{-2}\int_0^\theta\big(\varphi(r\cos\phi)f(\omega r\sin\phi)+g(r\cos\phi)\big)\sin\phi d\phi.$$
Under this transformation, (\ref{appl-syst4}) is transformed into
\begin{equation}\label{appl-syst5}
\begin{cases}
\frac{d\lambda}{d\theta}=-\omega^{-2}p(t)\sin\theta+O_{3m+2}(\lambda^{-1}),\\
\frac{dt}{d\theta}=\omega^{-1}-\omega^{-3}\lambda^{-1}\big(\varphi(\lambda\cos\theta)f(\omega \lambda\sin\theta)+g(\lambda\cos\theta)\big)\cos\theta+\omega^{-3}\lambda^{-1}p(t)\cos\theta+O_{3m+3}(\lambda^{-2}).
\end{cases}
\end{equation}

Introduce a transformation $$\lambda=\lambda,\quad \tau=t+S_2(\theta,\lambda),$$
where $$S_2(\theta,\lambda)=\omega^{-3}\lambda^{-1}\int_0^\theta\big(\varphi(\lambda\cos\phi)f(\omega \lambda\sin\phi)\cos\phi-\lambda J_1(\lambda)\big)+\big(g(\lambda\cos\phi)\cos\phi-\lambda J_2(\lambda)\big) d\phi,$$
with $$J_1(\lambda)=\frac{1}{2\pi\lambda}\int_0^{2\pi}\varphi(\lambda\cos\phi)f(\omega \lambda\sin\phi)\cos\phi d\phi,$$
$$J_2(\lambda)=\frac{1}{2\pi\lambda}\int_0^{2\pi}g(\lambda\cos\phi)\cos\phi d\phi.$$
By this transformation, (\ref{appl-syst5}) is transformed into
\begin{equation}\label{appl-syst6}
\begin{cases}
\frac{d\lambda}{d\theta}=-\omega^{-2}p(\tau)\sin\theta+O_{3m+2}(\lambda^{-1})\\
\frac{d\tau}{d\theta}=\omega^{-1}-\omega^{-3}\big(J_1(\lambda)+J_2(\lambda)\big)+\omega^{-3}\lambda^{-1}p(\tau)\cos\theta+O_{3m+2}(\lambda^{-2}).
\end{cases}
\end{equation}
Furthermore, we can find a transformation $(\lambda,\tau)\rightarrow(\lambda,\varsigma)$,
where \begin{equation}\label{change2}
\varsigma=\tau+\lambda^{-1}S_3(\theta,\tau),
\end{equation}
and $S_3(\theta,\tau)$ is determined by solving equation
\begin{equation}\label{homol3}
\omega^{-3}p(\tau)\cos\theta+\frac{\partial S_3}{\partial\theta}+\omega^{-1}\frac{\partial S_3}{\partial\tau}=0.
\end{equation}

Write $p(\tau)$ in the form of Fourier series
\begin{equation}\label{series1}p(\tau)=\sum_{k\in\mathbb{Z}^m}p_ke^{i\langle k,\mu\rangle\tau},
\end{equation}
and $p_{-k}=\overline{p_k}$.

In addition,
 \begin{equation}\label{series4}\cos\theta=\frac{1}{2}(e^{i\theta}+e^{-i\theta}),
\end{equation}
and denote
\begin{equation}\label{series5}
S_3(\theta,\tau)=e^{i\theta}\sum_{k\in\mathbb{Z}^m}\chi_k^+e^{i\langle k,\mu\rangle\tau}+e^{-i\theta}\sum_{k\in\mathbb{Z}^m}\chi_k^-e^{i\langle k,\mu\rangle\tau},
\end{equation}
Substituting (\ref{series1}), (\ref{series4}) and (\ref{series5}) into (\ref{homol3}), one has that
$$S_3(\theta,\tau)=e^{i\theta}\sum_{k\in\mathbb{Z}^m}\frac{i\omega^{-3}}{2\big(\langle k,\mu\rangle\omega^{-1}+1\big)}p_ke^{i\langle k,\mu\rangle\tau}+e^{-i\theta}\sum_{k\in\mathbb{Z}^m}\frac{i\omega^{-3}}{2\big(\langle k,\mu\rangle\omega^{-1}-1\big)}p_ke^{i\langle k,\mu\rangle\tau}.$$
 Different from the real analytic case, we need to prove the series $S_3(\theta,\tau)$ is convergent.

 From $p\in\mathcal{C}^{4m+5}$ and (\ref{diophantine}),  we note that $|p_k|\leq\|p\|_{\mathcal{C}^{4m+5}}|k|^{-4m-5}$ and $$\left|\frac{\omega^{-2}}{2\big(\langle k,\mu\rangle\omega^{-1}+1\big)}p_k\right|\leq C\|p\|_{\mathcal{C}^{4m+5}}|k|^{-4m-5+\sigma}.$$ By $m<\sigma<m+1$, it yields that the series converges to $S_3(\theta,\tau)$. Moreover, $S_3$ is $\mathcal{C}^\infty$, $2\pi$-periodic in $\theta$ and $\mathcal{C}^{3m+3}$, quasi-periodic in $\tau$ with the frequency $\mu$. $S_3\in O_{3m+3}(1)$, $\frac{\partial S_3}{\partial \tau}\in O_{3m+2}(1)$. Therefore, we derive the system
 \begin{equation}\label{appl-syst7}
\begin{cases}
\frac{d\lambda}{d\theta}=-\omega^{-2}p(\varsigma)\sin\theta+O_{3m+2}(\lambda^{-1})\\
\frac{d\varsigma}{d\theta}=\omega^{-1}-\omega^{-3}\big(J_1(\lambda)+J_2(\lambda)\big)+O_{3m+2}(\lambda^{-2}).
\end{cases}
\end{equation}
Thus by this transformation, we eliminate $\omega^{-3}\lambda^{-1}p(\tau)\cos\theta$ item in the second equation of (\ref{appl-syst6}).

Since $$\lim_{\lambda\rightarrow+\infty}\lambda^{k+1}J_1^{(k)}(\lambda)=(-1)^kk!\frac{1}{\pi}\big(\varphi(+\infty)-\varphi(-\infty)\big)f(+\infty),$$ and
$$\lim_{\lambda\rightarrow+\infty}\lambda^{k+1}J_2^{(k)}(\lambda)=(-1)^kk!\frac{1}{\pi}(g(+\infty)-g(-\infty)), \quad0\leq k\leq 3m+3,$$
 (\ref{appl-syst6}) can be rewritten as
\begin{equation}\label{appl-syst7}
\begin{cases}
\frac{d\lambda}{d\theta}=-\omega^{-2}p(\varsigma)\sin\theta+O_{3m+2}(\lambda^{-1})\\
\frac{d\varsigma}{d\theta}=\omega^{-1}-\frac{\omega^{-3}}{\pi}\lambda^{-1}\Big(\big(\varphi(+\infty)-\varphi(-\infty)\big)f(+\infty)+\big(g(+\infty)-g(-\infty)\big)\Big)+o_{3m+2}(\lambda^{-1}).
\end{cases}
\end{equation}
We also recalled (\ref{appl-syst7}) is reversible with respect to $(\lambda,\varsigma)\mapsto(\lambda,-\varsigma)$.

Denote $\lambda=\rho^{-1},$ then (\ref{appl-syst7}) can be rewritten as

\begin{equation*}
\begin{cases}
\frac{d\rho}{d\theta}=\omega^{-2}\rho^2p(\varsigma)\sin\theta+o_{3m+2}(\rho^2),\\
\frac{d\varsigma}{d\theta}=\omega^{-1}-\frac{\omega^{-3}}{\pi}\rho\Big(\big(\varphi(+\infty)-\varphi(-\infty)\big)f(+\infty)+\big(g(+\infty)-g(-\infty)\big)\Big)+o_{3m+2}(\rho).
\end{cases}
\end{equation*}
Therefore, we derive the corresponding Poincar\'{e} mapping with the form
 \begin{equation}\label{appl-poin}
\begin{cases}
\rho_1=\rho_0+\rho_0^2 l(\varsigma_0)+o_{3m+2}(\rho_0^2),\\
\varsigma_1=\varsigma_0+\gamma_0+\gamma_1\rho_0+o_{3m+2}(\rho_0),
\end{cases}
\end{equation}
where $\gamma_0=2\pi\omega^{-1}$, $l(\varsigma_0)=\omega^{-2}\int_0^{2\pi} p(\varsigma_0+\omega^{-1}\theta)\sin\theta d\theta$
and
 \begin{equation}\label{twistcondition}
 \gamma_1=-2\omega^{-3}\big((\varphi(+\infty)-\varphi(-\infty))f(+\infty)+(g(+\infty)-g(-\infty))\big).
 \end{equation}

If $(\varphi(+\infty)-\varphi(-\infty))f(+\infty)+(g(+\infty)-g(-\infty))\neq0$,  by Remark \ref{remark}, the mapping (\ref{appl-poin}) has many invariant curves tending to $\rho_0=0$, which means the invariant curves of the Poincar\'{e} map of (\ref{appl-equa}) tend to infinity. Thus for equation (\ref{appl-equa}), the existence of quasi-periodic solutions is got. Moreover for the initial value lying between two invariant curves, the solution is globally bounded. As the invariant curves tend to infinity, all solutions of (\ref{appl-equa}) are bounded, therefore we finish the proof of Theorem \ref{appl-theo}.

\section*{Acknowledgment}

Y. Zhuang was partially supported by the NSFC (grant no. 11971059) and the Fundamental Research Funds for the Central Universities (grant no. 202261096).
D. Piao was supported in part by the NSFC (grant no. 11971059).  Y. Niu was supported by the NSFC (grant no. 12201587) , Fundamental Research Funds for the Central Universities (grant no. 202213040) and the Shandong Provincial Natural Science Foundation, China (grant no. A010704).



\end{spacing}

\end{document}